\documentclass[11pt]{amsart}
\usepackage[latin1]{inputenc}
\usepackage[T1]{fontenc}
\usepackage[english,francais]{babel}
\usepackage{ae, aecompl}
\usepackage{amsmath,amsfonts,amssymb}
\usepackage{amsthm}
\usepackage{mathrsfs,array,graphicx}
\usepackage{stmaryrd}
\usepackage{wasysym}
\usepackage[all,ps]{xy}
\usepackage{calc}
\usepackage{enumerate}
\usepackage{ulem}

\usepackage{fancyhdr}
\addto\captionsfrench{}

\newcommand{\ugu}{\stackrel{\mathrm{d\acute{e}f}}{=}}
\newcommand{\unif}{\rho(\varpi_F)}

\newcommand{\Q}{\mathbb{Q}_p}

\newcommand{\Z}{\mathbb{Z}_p}

\newcommand{\OF}{\mathcal{O}_F}

\newcommand{\uni}{\varpi_F}
\newcommand{\C}{C}

\newcommand{\Cr}{C^{r}}

\newcommand{\into}{\hookrightarrow}

\newcommand{\GL}{\mathrm{GL}_2}

\newcommand{\Gal}{\mathrm{Gal}}

\newtheoremstyle{Proposition}{11pt}{11pt}{\itshape}{}{\bfseries}{.}{.5em }{}
\theoremstyle{Proposition}

\newcounter{Proposition}
\newtheorem{theo}[Proposition]{Théorème}
\newtheorem{prop}[Proposition]{Proposition}

\newtheorem{defin}[Proposition]{Définition}
\newtheorem{exam}[Proposition]{Exemple}

\newtheorem{rem}[Proposition]{Remarque}
\newtheorem{cor}[Proposition]{Corollaire}
\newtheorem{lemma}[Proposition]{Lemme}

\numberwithin{equation}{section}
\numberwithin{Proposition}{section}

\input xy
\xyoption{all}
\usepackage{amscd}

\author{Marco De Ieso}
\address{Bâtiment 430, Université Paris-Sud, 91405, Orsay Cedex, France}
\email{Marco.DeIeso@math.u-psud.fr}

\title{ESPACES DE FONCTIONS DE CLASSE $C^r$ SUR $\OF$}

\begin{document}
\date{}
\pagestyle{fancy}

\







\begin{otherlanguage}{english}
\begin{abstract}
In this paper we introduce a class of Banach spaces of functions of class $C^r$ (where $r$ is a positive real number) and the associated dual spaces of distributions of order $r$, which turn out to be useful in $p$-adic Langlands theory (\cite{marco7}). We construct a Banach basis for these spaces and we give a criterion for telling when a linear form on a space of locally $\Q$-polynomial functions extends to a distribution of order $r$. This generalises the classical results of Amice-Vélu and Vishik (\cite{amivel}, \cite{vis}).
\vskip.5cm
  \par\noindent   \textsc{R{\'e}sum{\'e}.} Dans cet article nous introduisons certains espaces de Banach de fonctions de classe $C^r$, où $r$ est un réel positif, et leurs duaux des distributions d'ordre $r$, qui se revèlent utiles en théorie de Langlands $p$-adique (\cite{marco7}).  Nous construisons une base de Banach de ces espaces et nous donnons un critère pour qu'une forme linéaire sur un  espace de fonctions localement $\Q$-polynomiales s'étende en une distribution d'ordre $r$, ce qui généralise des résultats classiques d'Amice-Vélu et Vishik (\cite{amivel}, \cite{vis}).
\end{abstract}
\end{otherlanguage}

\maketitle

\tableofcontents

\section{Introduction} 

Soit $p$ un nombre premier. La correspondance de Langlands $p$-adique pour $\GL(\Q)$, commencée sous l'impulsion de Breuil (\cite{breuilab}, \cite{breuila}) et établie par  Colmez \cite{colmez2} et  Pa\v{s}k\={u}nas \cite{pask} à la suite des travaux de Colmez \cite{colmez3} et Berger-Breuil \cite{bb}, est une bijection entre certaines représentations de dimension $2$ de $\Gal(\overline{\mathbb{Q}}_p/\Q)$ et certaines représentations de $\GL(\Q)$. En particulier, dans \cite{bb}, une étape importante dans l'établissement de cette correspondance pour la classe des représentations galoisiennes de dimension $2$, absolument irréductibles et devenant cristallines sur une extension abélienne de $\Q$, est la description explicite du complétés unitaire universel de certaines induites paraboliques localement algébriques (\cite[Théorème 4.3.1]{bb}). C'est là que l'analyse $p$-adique, et plus précisement l'espace de Banach des fonctions de classe $C^r$ sur $\Z$, où $r$ désigne un nombre réel positif et $\Z$ l'anneau des entiers $p$-adiques, intervient de manière cruciale. Une question naturelle est la généralisation de cette étape pour certaines représentations localement $\Q$-analytiques de $\GL(F)$, où $F$ est une extension finie de $\Q$ (\cite{marco7}).

Pour cela nous avons été amené à introduire une nouvelle notion de fonction de classe $C^r$ sur $\OF$ que nous développons dans cet article.

Rappelons que pour $\OF = \Z$ trois constructions \textit{a priori} différentes de l'espace des fonctions de classe $C^r$ à valeurs dans $E$, où $E$ désigne une extension finie de $\Q$, ont été introduites par Barsky et Schikhof dans le cas où $r$ est un entier positif  (\cite{ba}, \cite{sk2}), et ont été généralisées pour $r$ un réel positif quelconque grâce aux travaux de Berger-Breuil, Colmez et Nagel (\cite{bb}, \cite{enno1}, \cite{colmez}). On sait maintenant que ces trois constructions donnent des notions équivalentes (\cite{sk}, \cite{enno1}). La définition que l'on va généraliser dans le cas où $F$ est une extension finie de $\Q$ est celle de (\cite{sk}, \cite{colmez}) que nous rappelons maintenant. Fixons un réel positif $r$ et notons $[r]$ sa partie entière. Une fonction $f\colon \mathbb{Z}_p \to E$ est dite de classe $C^r$ si $f(x+y)$ a un développement limité à l'ordre $[r]$ en tout $x$, et si le reste est $o(|y|^r)$ uniformément (en $x$) sur tout compact.

Le premier résultat principal de cet article est la construction d'une base de Banach de l'espace  $C^r(\OF,E)$ des fonctions de classe $C^r$ sur $\OF$, qui consiste d'une famille dénombrable de fonctions localement $\Q$-polynomiales. Le deuxième, qui porte sur son dual topologique, donne une condition nécessaire et suffisante pour qu'une forme linéaire définie sur un espace de fonctions localement $\Q$-polynomiales convenable s'étende en une distribution d'ordre $r$, ce qui généralise des résultats classiques d'Amice-Vélu et Vishik (\cite{amivel}, \cite{vis}).

\subsection{Notations} Soit $p$ un nombre premier. On fixe une clôture algébrique $\overline{\mathbb{Q}}_p$ de $\Q$ et une extension finie $F$ de $\Q$ contenue dans $\overline{\mathbb{Q}}_p$. On désignera toujours par $E$ une extension finie de $\Q$ qui vérifie:
\[
|S| = [F:\Q],
\]
où $S \ugu \mathrm{Hom}_{alg}(F,E)$. 

En général, si $L$ désigne $F$ ou $E$, on note $\mathcal{O}_L$ son anneau d'entiers, $\varpi_L$ une uniformisante de $\mathcal{O}_L$ et $k_L = \mathcal{O}_L/(\varpi_L)$ son corps résiduel. On note $f = [k_F : \mathbb{F}_p]$, $q= p^f$ et $e$ l'indice de ramification de $F$ sur $\Q$, de sorte que $[F:\Q] = ef$ et $k_F \simeq \mathbb{F}_q$.  

La valuation $p$-adique $val_F$  sur $\overline{\mathbb{Q}}_p$ est normalisée par $val_F(p) = [F:\Q]$ et on pose $|x| = p^{- val_F(x)}$  si $x \in \overline{\mathbb{Q}}_p$.

Si $n \in \mathbb{Z}_{\geq 0}$ et $* \in \{<,\leq,>, \geq, = \}$ notons: 
\[
I_{*n} = \Big\{(i_{\sigma})_{\sigma \in S} \in \mathbb{Z}_{\geq 0}^{|S|}, \, \sum_{\sigma \in S}i_{\sigma} * n  \Big\}.
\]
Si $\underline{n} = (n_{\sigma})_{\sigma \in S}, \underline{m} = (m_{\sigma})_{\sigma \in S}$ sont des $|S|$-uplets d'entiers positifs ou nuls posons: 
\begin{itemize}
\item[(i)] $\underline{n}! = \prod_{\sigma \in S}n_{\sigma}!$;
\item[(ii)] $|\underline{n}| = \sum_{\sigma \in S}n_{\sigma}$;
\item[(iii)] $\underline{n}-{\underline{m}} = (n_{\sigma}-m_{\sigma})_{\sigma \in S}$; 
\item[(iv)] $\underline{n}\leqslant{\underline{m}}$ si $n_{\sigma}\leq m_{\sigma}$ pour tout $\sigma \in S$; 
\item[(v)] $\binom{\underline{n}}{\underline{m}} = \frac{\underline{n}!}{\underline{m}!(\underline{n}-\underline{m})!}$.
\end{itemize} 
Si $\underline{n}= (n_{\sigma})_{\sigma \in S} \in \mathbb{Z}_{\geq 0}^{|S|}$ et $z \in \OF$ on pose $z^{\underline{n}} = \prod_{\sigma \in S}\sigma(z)^{n_{\sigma}}$.

Une norme $p$-adique sur un $E$-espace vectoriel $V$ est une fonction $\|\cdot\| \colon V \to \mathbb{R}_{\geq 0}$ telle que: 
\begin{itemize}
\item[(i)] $\|v+w\| \leq \sup(\|v\|,\|w\|)$ pour tout $v,w \in V$;
\item[(ii)] $\|\lambda v\| \leq |\lambda| \|v\|$ pour tout $\lambda \in E$, $v\in V$;
\item[(iii)] $\|v \| = 0$ si et seulement si $v = 0$.
\end{itemize}
Un espace de Banach $p$-adique sur $E$ est un $E$-espace vectoriel topologique complet dont la topologie provient d'une norme $p$-adique. Dans ce texte tous les espaces de Banach sont $p$-adique et tels que $\|B \| \subseteq |E|$.


Si $V$ est un $E$-espace vectoriel topologique, on note $V^{\vee}$ son dual topologique.

\subsection{Énoncé des résultats}

Dans la Section $2$ nous introduisons d'abord la notion de fonction de classe $C^r$ sur $\OF$, où $r$ est un réel positif. On dit qu'une fonction  $f\colon \OF \to E$ est de classe $\Cr$ s'il existe une famille de fonctions bornées $D_{\underline{i}}f \colon \OF \to E$, pour $\underline{i} \in I_{\leq [r]}$, telles que $f(x+y)$ a un développement limité à l'ordre $[r]$ en tout $x$ et si le reste est $o(|y|^r)$ uniformément (en $x$) sur tout compact.   Montrer l'unicité de cette famille de fonctions requiert une estimation technique et donnée en appendice, sur le maximum des valeurs absolues des coefficients dominants d'une fonction $\Q$-algébrique (Proposition \ref{ultratec1}). Dans une deuxième partie on étudie quelques propriétés de l'espace $C^r(\OF,E)$: nous montrons que c'est une algèbre de Banach sur $E$ et que pour tout $\underline{i} \in I_{\leq [r]}$ l'opérateur  $D_{\underline{i}}$ définit une application  continue de $C^r$ dans $C^{r-|\underline{i}|}$. Ensuite nous donnons une condition suffisante sur une fonction $h\colon \OF \to \OF$ pour que $f\circ h$ soit de classe $C^r$ et nous montrons que si l'on fixe une telle fonction alors l'application qui associe à toute $f \in C^r(\OF,E)$ la fonction $f\circ h$ est continue.  
\medbreak

Soit $J$ une partie de $S$ et $(d_{\sigma})_{\sigma \in S\setminus J}$ un $|S\backslash J|$-uplet d'entiers positifs ou nuls. Nous considérons dans la section $3$ le sous-espace $\mathcal{F}(\OF,J,(d_{\sigma})_{\sigma \in S\setminus J})$ de l'espace des fonctions localement $\Q$-analytiques qui sont localement analytiques selon les plongements $\sigma$ dans $J$ et localement polynomiales de degré au plus $d_{\sigma}$ selon les plongements $\sigma$ dans $S\backslash J$. Nous démontrons que cet espace s'injecte de façon continue dans l'espace des fonctions de classe $C^r$ (Corollaire \ref{inicont}) et nous en décrivons l'adhérence dans $C^r(\OF,E)$ en utilisant les opérateurs de dérivation $D_{\underline{i}}$ (Corollaire \ref{densità2}). On note $C^r(\OF,J',(d_{\sigma})_{\sigma \in S\setminus J'})$ cette adhérence, où $J'$ désigne le sous-ensemble de $S$ défini par:
\[
J' = J \coprod \{\sigma \in S\backslash J,\, d_{\sigma}+1 > r\}.
\]   
Une deuxième partie de cette section est consacrée à la construction d'une base de Banach (qui dépend de $r$) de l'espace $C^r(\OF,E)$. Dans le cas $\OF = \Z$ mentionnons que cette base coïncide avec celle construite par Van der Put pour l'espace des fonctions continues  sur $\Z$ (ou ce qui revient au même pour $C^0(\Z,E)$) et généralisée par Colmez pour $C^r(\Z,E)$ avec $r$ un réel quelconque (\cite{vander}, \cite{colmez}). Signalons que pour l'espace des fonctions continues sur $\OF$ cette base a déjà été construite par De Shalit (\cite[§2]{sha}). Donnons un aperçu de notre construction. Fixons un plongement $\rho \colon F \into E$. Posons $A_0 = \{0\}$, choisissons pour tout $h \in \mathbb{Z}_{>0}$ un système de représentants $A_h \subset \OF$ des classes de  $\OF/\varpi_F^h \OF$ de sorte que $A_h \supset A_{h-1}$ et notons $A = \coprod_{h \geq 1} A_h \backslash A_{h-1}$. Pour tout $a \in A$ on définit $l(a)$ comme le plus petit entier $n_0$ tel que $a \in A_{n_0}$.   Si $a \in A$ et $\underline{i} \in I_{\leq [r]}$, on note $e_{a,\underline{i},r}$ la fonction définie par:
\[
z\mapsto e_{a,\underline{i},r} (z) = \unif^{[l(a)r]} \mathbf{1}_{a+\varpi_F^{l(a)} \OF}(z) \Big(\frac{z-a}{\varpi_F^{l(a)}}\Big)^{\underline{i}}.
\]

\begin{theo}\label{princ}
La famille des $e_{a,\underline{i},r}$, pour $a \in A$ et $\underline{i} \in I_{\leq [r]}$, forme une base de Banach de $C^r(\OF,E)$.
\end{theo}
Pour cela on généralise la preuve de \cite[I.5.14]{colmez}. Plus précisement, un premier ingrédient que nous utilisons est une estimation de la norme $C^r$ des $e_{a,\underline{i},r}$ (Lemme \ref{lemval}). Un deuxième est la construction explicite, pour toute fonction $f\in C^r(\OF,E)$, d'une suite de fonctions $f_h$ (qui dépendent de $f$) telles que $f_h$ tend vers $f$ dans $C^r(\OF,E)$ quand $h$ tend vers $+\infty$ (Proposition \ref{densità}). Nous terminons cette Section en décrivant la sous-famille de la famille des $e_{a,\underline{i},r}$ qui va être une base de Banach pour l'espace $C^r(\OF,J',(d_{\sigma})_{\sigma \in S\setminus J'})$. 
\medbreak

Dans la Section $4$ on s'occupe des duaux topologiques des espaces considérés précédemment. Si $N\in \mathbb{Z}_{\geq 0}$, on note $\mathcal{F}^N(\OF,J, (d_{\sigma})_{\sigma \in S \setminus J})$ l'espace des fonctions localement $\Q$-algébriques sur $\OF$ de degré au plus $N$ dans $\mathcal{F}(\OF,J, (d_{\sigma})_{\sigma \in S \setminus J})$ et $\mathcal{F}^N(\OF,J, (d_{\sigma})_{\sigma \in S \setminus J})^{\vee}$ l'ensemble des formes linéaires sur $\mathcal{F}^N(\OF,J, (d_{\sigma})_{\sigma \in S \setminus J})$. Le résultat principal (Théorème \ref{velu}) donne une condition nécessaire et suffisante pour qu'une forme linéaire sur $\mathcal{F}^N(\OF,J, (d_{\sigma})_{\sigma \in S \setminus J})$ s'étende en une forme linéaire continue sur $C^r(\OF,J', (d_{\sigma})_{\sigma \in S \setminus J'})$. Cela généralise un résultat dû à Amice-Vélu et Vishik (\cite{amivel}, \cite{vis}). 
\begin{theo}\label{principes}
(i) Soit $\mu \in C^r(\OF,J', (d_{\sigma})_{\sigma \in S \setminus J'})^{\vee}$. Il existe une constante $C_{\mu} \in \mathbb{R}_{\geq 0}$ telle que pour tout $a \in \OF$, tout $n \in \mathbb{Z}_{\geq 0}$ et tout  $\underline{i} \in \mathbb{Z}_{\geq 0}^{|S|}$, $i_{\sigma} \leq d_{\sigma}$ pour tout $\sigma \in S\backslash J'$ on ait: 
\begin{align*}
\Big|\int_{a+\uni^n\OF} \Big(\frac{z-a}{\varpi_F^n} \Big)^{\underline{i}} \mu(z)\Big| \leq C_{\mu} \, q^{nr}. 
\end{align*}

(ii) Soit $N\geq [r]$ et $\mu \in  \mathcal{F}^N(\OF,J, (d_{\sigma})_{\sigma \in S \setminus J})^{\vee}$. Supposons qu'il existe une constante $C_{\mu} \in \mathbb{R}_{\geq 0}$ telle que, pour tout $a \in \OF$, tout $n \in \mathbb{Z}_{\geq 0}$ et tout $\underline{i} \in I_{\leq N}$, $i_{\sigma} \leq d_{\sigma}$ pour tout $\sigma \in S\backslash J$ on ait:    
\begin{align*}
\Big|\int_{a+\uni^n\OF}  \Big(\frac{z-a}{\varpi_F^n} \Big)^{\underline{i}} \mu(z) \Big| \leq C_{\mu} \, q^{nr}.
\end{align*}
Alors $\mu$ se prolonge de manière unique en une forme linéaire continue sur $C^r(\OF,J', (d_{\sigma})_{\sigma \in S \setminus J'})$. 
\end{theo}

Le (i) est conséquence d'une estimation de la norme $C^r$ d'une fonction localement $\Q$-analytique (Proposition \ref{grossaz}). Le (ii) est plus subtil et est basé sur deux ingrédients. Le premier est la base de Banach de l'espace $C^r(\OF,J', (d_{\sigma})_{\sigma \in S \setminus J'})$ que l'on construit dans la Section 3. Le deuxième est un résultat de densité (Corollaire \ref{densità2}).    
\medbreak

La Section 5 est consacrée à l'étude d'une autre notion de fonction de classe $C^r$ sur $\OF$ qui s'appuie sur le fait que $\OF$ est un $\Z$-module libre de rang fini. On démontre que cette deuxième notion n'est pas équivalente à la première.
\medbreak

{\scshape{Remerciements}}. Je remercie chalheuresement mon directeur de thèse Christophe Breuil pour ses conseils, pour ses très nombreuses remarques et pour avoir suivi attentivement l'évolution de ce travail. Je remercie Benjamin Schraen pour avoir répondu à mes questions et pour avoir  lu avec intérêt une version préliminaire de ce travail. Je remercie Arno Kret et Enno Nagel pour des discussions utiles.

\addtocontents{toc}{\protect\setcounter{tocdepth}{2}}

\section{Fonctions de classe $C^r$ sur $\OF$}

\subsection{Définition}\label{classe}

Soit $r \in \mathbb{R}_{\geq 0}$. Notons $[r]$ sa partie entière. On dispose de l'espace de Banach de fonctions de classe $C^r$ sur $\Z$ à valeurs dans $E$. Rappelons (\cite{colmez}) que $f\colon \mathbb{Z}_p \to E$ est de classe $C^r$ si $f(x+y)$ a un développement limité à l'ordre $[r]$ en tout $x$, et si le reste est $o(|y|^r)$ uniformément (en $x$) sur tout compact. Dans cette section on va construire un espace de fonctions de $\OF$ dans $E$ qui généralise cette idée. La définition que l'on donne ne dépend pas du choix d'une $\Z$-base de $\OF$: elle va dépendre juste des plongements de $\OF$ dans $E$.


\begin{defin}\label{definizione}
On dit que $f\colon \OF \to E$ est de classe $\Cr$ sur $\OF$ s'il existe des fonctions bornées $D_{\underline{i}}f \colon \OF \to E$, pour $\underline{i} \in I_{\leq [r]}$, telles que, si l'on définit $\varepsilon_{f,[r]}\colon \OF \times \OF \to E$ par: 
\begin{align*}
\varepsilon_{f,[r]}(x,y) = f(x+y) -
 \sum_{\underline{i} \in I_{\leq [r]}} D_{\underline{i}}f(x) \frac{y^{\underline{i}}}{\underline{i}!}
\end{align*}
et pour tout $h \in \mathbb{Z}_{\geq 0}$
\begin{align*}
C_{f,r}(h) = \sup_{x \in \OF, y \in \varpi_F^h \OF} |\varepsilon_{f,[r]}(x,y)|  q^{rh}
\end{align*}
alors $C_{f,r}(h)$ tend vers $0$ quand $h$ tend vers $+\infty$. 
\end{defin}

\begin{rem}\label{unifcont}
{\rm (i) Soit $f$ une fonction de classe $C^r$ sur $\OF$. Alors $\varepsilon_{f,[r]}(x,0)=0$ pour tout $x \in \OF$ car, pour tout $h$ on a: 
\[
\sup_{x \in \OF}|\varepsilon_{f,[r]}(x,0)| \leq \sup_{x \in \OF}|\varepsilon_{f,[r]}(x,0)|  q^{rh} \leq C_{f,r}(h)
\] 
et $C_{f,r}(h)$ tend vers $0$ quand $h$ tend vers $+\infty$ par hypothèse. Cela implique  $D_{\underline{0}}f(x) = f(x)$ pour tout $x \in \OF$. Comme $C_{f,r}(h)$ tend vers $0$ quand $h \to +\infty$, \textit{a fortiori} on a:
\begin{align}\label{revisione4}
\sup_{x \in \OF, y \in \varpi_F^h \OF} |\varepsilon_{f,[r]}(x,y)| \to 0 \ \mbox{quand} \ h\to +\infty.
\end{align}
De plus: 
\begin{align}\label{revisione3}
\sup_{x \in \OF, y \in \varpi_F^h \OF}
\Bigl|\sum_{\substack{
 \underline{i} \in I_{\leq [r]} \\
|\underline{i}| \geq 1  
}}  D_{\underline{i}}f(x) \frac{y^{\underline{i}}}{\underline{i}!}\Bigr| \to 0 \ \mbox{quand} \ h \to +\infty
\end{align}   
car par définition $D_{\underline{i}}f$ est une fonction bornée sur $\OF$ pour tout $\underline{i} \in I_{\leq [r]}$  et $|y^{\underline{i}}| \leq q^{-h}$ pour tout $y \in \varpi_F^h \OF$ et tout $\underline{i} \in I_{\leq [r]}, |\underline{i}| \geq 1$.
On peut réécrire $\varepsilon_{f,[r]}(x,y)$ sous la forme
\[
\varepsilon_{f,[r]}(x,y) = \big(f(x+y)-f(x)\big) - \sum_{\substack{
 \underline{i} \in I_{\leq [r]} \\
|\underline{i}| \geq 1  
}}  D_{\underline{i}}f(x) \frac{y^{\underline{i}}}{\underline{i}!}
\]
et donc \eqref{revisione4} et \eqref{revisione3} impliquent
\[
\sup_{x \in \OF, y \in \varpi_F^h \OF}
|\ f(x+y) - f(x) | \to 0 \ \mbox{quand} \ h \to +\infty
\]
c'est-à-dire $f$ est une fonction continue. En particulier une fonction de classe $C^0$ est une fonction continue. 

(ii) Soit $f$ une fonction continue sur $\OF$. Alors $f$ est de classe $C^0$ sur $\OF$. En effet, par continuité de $f$ on a:
\[
\sup_{x \in \OF, y \in \varpi_F^h \OF}
|\ f(x+y) - f(x) | \to 0 \ \mbox{quand} \ h \to +\infty
\]
et donc si l'on pose: 
\[ 
\forall x \in \OF, \quad D_{\underline{0}}f(x) = f(x) 
\]
on voit que $D_{\underline{0}}f$ vérifie la Définition \ref{definizione} et cela permet de conclure.}
\end{rem}

\begin{lemma}\label{cl}
Soit $f\colon \OF \to E$ une fonction de classe $C^r$. Alors $f$ est de classe $C^l$ pour tout $l \in \mathbb{R}_{\geq 0}$, $l\leq r$.
\begin{proof}
Il suffit de prouver que $f$ est de classe $C^l$ pour tout $l \in \mathbb{Z}_{\geq 0}$ et $l \leq [r]$. La démonstration se fait par récurrence sur $l\leq [r]$. Supposons donc que $f$ soit de classe $C^l$ pour un $l\in \mathbb{Z}_{> 0}$, $l \leq [r]$. Montrons qu'elle est de classe $C^{l-1}$. Comme $C_{f,l}(h)$ tend vers $0$ quand $h$ tend vers $+\infty$, alors \textit{a fortiori} on a
\begin{align}\label{dec1}
\sup_{x \in \OF, y \in \varpi_F^h \OF} |\varepsilon_{f,[l]}(x,y)|  q^{(l-1)h} \to 0 \quad \mbox{quand} \quad h \to +\infty.
\end{align}
D'autre part  
\begin{align}\label{dec2}
\sup_{x \in \OF, y \in \varpi_F^h \OF}
\Bigl|\sum_{\underline{i} \in I_{=l}}  D_{\underline{i}}f(x)  \frac{y^{\underline{i}}}{\underline{i}!}\Bigr|  q^{h(l-1)} \to 0 \ \mbox{quand} \ h \to +\infty.
\end{align}  
En effet, la fonction $D_{\underline{i}}f$ est bornée sur $\OF$ pour tout $\underline{i} \in I_{=l}$ par définition et $|y^{\underline{i}}| \leq q^{-hl}$ pour tout $y \in \varpi_F^h \OF$, $\underline{i} \in I_{=l}$. Donc si l'on pose
\[
C = \sup_{\underline{i} \in I_{=l}} \sup_{x \in \OF} \big|D_{\underline{i}}f(x) \big|
\]
on obtient:
\[
\sup_{x \in \OF, y \in \varpi_F^h \OF}
\Bigl|\sum_{\underline{i} \in I_{=l}}  D_{\underline{i}}f(x)  \frac{y^{\underline{i}}}{\underline{i}!}\Bigr|  q^{h(l-1)} \leq C q^{-h}
\]
et $C q^{-h} \to 0$ quand $h \to +\infty$. On peut réécrire $\varepsilon_{f,[l]}(x,y)$ sous la forme: 
\[
\varepsilon_{f,[l]}(x,y) = \varepsilon_{f,[l-1]}(x,y) - \sum_{\underline{i} \in I_{=l}}  D_{\underline{i}}f(x) \frac{y^{\underline{i}}}{\underline{i}!}.
\]
Les conditions \eqref{dec1} et \eqref{dec2} impliquent
\[
\sup_{x \in \OF, y \in \varpi_F^h \OF} |\varepsilon_{f,[l-1]}(x,y)|  q^{(l-1)h} \to 0 \quad \mbox{quand} \quad h \to +\infty,
\]
ce qui permet de conclure. 

\end{proof}
\end{lemma}

\begin{lemma}[Unicité]\label{unic}
Soit $f\colon \OF \to E$ une fonction de classe $C^r$. Alors il existe une unique famille de fonctions 
\[
\big\{ D_{\underline{i}}f \colon \OF \to E, \ \underline{i} \in I_{\leq [r]} \big\}
\]
qui vérifie la Définition \ref{definizione}.
\begin{proof}
Supposons que 
\begin{align*}
 \big\{F_{\underline{i}}\colon \OF \to E, \ \underline{i} \in I_{\leq [r]} \big\}  \quad \mbox{et} \quad
 \big\{G_{\underline{i}}\colon \OF \to E, \ \underline{i} \in I_{\leq [r]} \big\}
\end{align*}
soient deux familles de fonctions qui vérifient la Définition \ref{definizione}. On veut prouver que $F_{\underline{i}} = G_{\underline{i}}$ pour tout $\underline{i} \in I_{\leq [r]}$. La démonstration se fait par récurrence sur $|\underline{i}|$, $\underline{i} \in  I_{\leq [r]}$. D'après la Remarque \ref{unifcont} on a $F_{\underline{0}} = f= G_{\underline{0}}$. 

Soit $0 <  n < [r]$. Supposons $F_{\underline{i}} = G_{\underline{i}}$ pour tout $\underline{i} \in I_{\leq n-1}$. On veut montrer que $F_{\underline{i}} = G_{\underline{i}}$ pour tout $\underline{i} \in I_{=n}$. Par l'hypothèse de récurrence on déduit que pour tout $x,y \in \OF$ on a:
\[
f(x+y) -
 \sum_{\underline{i}\in I_{\leq n-1}} F_{\underline{i}}(x) \frac{y^{\underline{i}}}{\underline{i}!} = f(x+y) - \sum_{\underline{i}\in I_{\leq n-1}} G_{\underline{i}}(x) \frac{y^{\underline{i}}}{\underline{i}!}. 
\] 
Notons $\varphi_{n-1}$ la fonction définie pour tout $x,y \in \OF$ par:
\[
\varphi_{n-1}(x,y) = f(x+y) -
 \sum_{\underline{i}\in I_{\leq n-1}} F_{\underline{i}}(x) \frac{y^{\underline{i}}}{\underline{i}!}.
\]
D'après la preuve du Lemme \ref{cl} on sait que
\begin{align}\label{limit1}
\sup_{x \in \OF,y\in \varpi_F^h \OF}\Big|\varphi_{n-1}(x,y) -
 \sum_{\underline{i} \in I_{=n}} F_{\underline{i}}(x)  \frac{y^{\underline{i}}}{\underline{i}!}\Big|  q^{hn} \to 0 \ \mbox{quand} \ h \to +\infty
\end{align}
et 
\begin{align}\label{limit2}
\sup_{x \in \OF,y\in \varpi_F^h \OF}\Big|\varphi_{n-1}(x,y) -
 \sum_{\underline{i} \in I_{=n}} G_{\underline{i}}(x)  \frac{y^{\underline{i}}}{\underline{i}!}\Big|  q^{hn} \to 0 \ \mbox{quand} \ h \to +\infty.
\end{align}
Fixons $x \in \OF$. Supposons
\[
G_{\underline{i}}(x) = F_{\underline{i}}(x) + \alpha_{\underline{i}} \quad \alpha_{\underline{i}} \in E.
\]
Alors, par \eqref{limit2} on a:
\[
\sup_{y\in \varpi_F^h \OF}\Big|\Big(\varphi_{n-1}(x,y) -
\sum_{\underline{i} \in I_{=n}} F_{\underline{i}}(x)  \frac{y^{\underline{i}}}{\underline{i}!}\Big) + \sum_{\underline{i} \in I_{=n}}\alpha_{\underline{i}}\frac{y^{\underline{i}}}{\underline{i}!}\Big|  q^{hn} \to 0 \ \mbox{quand} \ h \to +\infty.
\]
Or, \eqref{limit1} implique
\begin{align}\label{disun}
\sup_{y\in \varpi_F^h \OF}\Big|\sum_{\underline{i} \in I_{=n}}\alpha_{\underline{i}}\frac{y^{\underline{i}}}{\underline{i}!}\Big|  q^{hn} \to 0 \ \mbox{quand} \ h \to +\infty.
\end{align}
Par la Proposition \ref{ultratec1} appliquée à $\sum_{\underline{i} \in I_{=n}}\alpha_{\underline{i}}\frac{y^{\underline{i}}}{\underline{i}!}$, il existe une constante $C \in \mathbb{R}_{\geq 1}$ et un $n_0 \in \mathbb{Z}_{\geq 0}$ tels que pour tout $h \geq n_0$ on a:
\[
\sup_{\underline{i} \in I_{=n}}|\alpha_{\underline{i}}| q^{-nh} \leq C \sup_{y \in \varpi_F^{h}\OF} \Big| \sum_{\underline{i} \in I_{=n}}\alpha_{\underline{i}}\frac{y^{\underline{i}}}{\underline{i}!} \Big| 
\]
et donc \eqref{disun} implique $|\alpha_{\underline{i}}|=0$ pour tout $\underline{i} \in I_{=n}$, d'où le résultat car ceci vaut pour tout $x \in \OF$.
\end{proof}
\end{lemma}

Notons $C^r(\OF,E)$ l'ensemble des fonctions de $\OF$ dans $E$ qui sont de classe $C^r$. On munit $C^r(\OF,E)$ de la norme $\|\cdot \|_{C^r}$ définie par:
\[
\|f \|_{C^r} = \sup \Big( \sup_{\underline{i} \in I_{\leq [r]}} \sup_{x \in \OF} \Bigl|\frac{D_{\underline{i}}f(x)}{\underline{i}!}\Bigr|, \sup_{x,y \in \OF} \frac{| \varepsilon_{f,[r]}(x,y)|}{|y|^{r}}   \Big)
\]
ce qui en fait un espace de Banach sur $E$.

\begin{exam}\label{esempio}
{\rm Soit $r \in \mathbb{R}_{\geq 0}$. Soit $\underline{m} = (m_{\sigma})_{\sigma \in S}$ un $|S|$-uplet d'entiers positifs. Considérons la fonction $f$ de $\OF$ dans $E$ définie par $z \mapsto z^{\underline{m}}$. Posons:
\begin{align}\label{esempp}
\frac{D_{\underline{i}}f(z)}{\underline{i}!} =
\left\{ \begin{array}{ll}
 \binom{\underline{m}}{\underline{i}} z^{\underline{m}-\underline{i}}   & \mbox{si} \ \underline{i} \leqslant \underline{m} \\
0 & \mbox{sinon.}
\end{array} \right.
\end{align}
Pour tout $x,y \in \OF$ on a:
\begin{align}\label{esempl}
\varepsilon_{f,[r]}(x,y) =
\left\{ \begin{array}{ll}
 \sum_{\substack{
 \underline{l} \leqslant \underline{m} \\
|\underline{l}| \geqslant [r] +1  
}} \binom{\underline{m}}{\underline{l}}y^{\underline{l}}x^{\underline{m}-\underline{l}}   & \mbox{si} \ r <|\underline{m}| \\
0 & \mbox{sinon.}
\end{array} \right.
\end{align}
En revenant à la Définition \ref{definizione} on voit que $f$ est bien de classe $C^r$. D'après \eqref{esempp} on a:
\[
\sup_{z \in \OF}\Big| \frac{D_{\underline{i}}f(z)}{\underline{i}!}\Big| \leq \sup_{z \in \OF} \Big|\binom{\underline{m}}{\underline{i}} z^{\underline{m}-\underline{i}}\Big| \leq 1  
\] 
et d'après \eqref{esempl} on a:
\[
\sup_{x,y \in \OF} \frac{| \varepsilon_{f,[r]}(x,y)|}{|y|^{r}} \leq \sup_{\substack{
 \underline{l} \leqslant \underline{m} \\
|\underline{l}| \geqslant [r] +1  
}} \sup_{x,y \in \OF}\frac{\big|\binom{\underline{m}}{\underline{l}}y^{\underline{l}}x^{\underline{m}-\underline{l}}\big|}{|y|^r} \leq 1.
\]
On en déduit que l'on a $\|f\|_{C^r} \leq 1$ et comme $\sup_{z\in \OF}|f(z)| = 1$ on a  $\|f\|_{C^r} = 1$.
}
\end{exam}

\subsection{Premières propriétés}\label{classe2}

Soit $r\in \mathbb{R}_{\geq 0}$ et $f$ une fonction de classe $C^r$ sur $\OF$. Alors, d'après le Lemme \ref{unic} il existe une unique famille de fonctions
\[
\big\{ D_{\underline{i}}f \colon \OF \to E, \ \underline{i}\in I_{\leq [r]} \big\}
\]
qui vérifie la Définition \ref{definizione}. Dans la suite on veut montrer que si l'on choisit une fonction de cette famille, disons $D_{\underline{i}}f$, alors elle est de classe $C^{r-|\underline{i}|}$. Nous aurons besoin du lemme technique suivant qui utilise une estimation, donnée en appendice, sur le maximum des valeurs absolues des coefficients dominants d'une fonction $\Q$-algébrique. 
\begin{lemma}\label{lemmatec}
Soit $N \in \mathbb{Z}_{\geq 0}$. Soit $a_{\underline{i}}$ pour $\underline{i} \in I_{\leq N}$ une famille d'éléments de $E$. Alors il existe une constante $C \in \mathbb{R}_{\geq 0}$ et un $n_0 \in \mathbb{Z}_{\geq 0}$ tels que pour tout $h \geq n_0$ on a:
\begin{align}\label{ineg}
\sup_{\underline{i} \in I_{\leq N}} |a_{\underline{i}}| q^{-h|\underline{i}|} \leq C \sup_{z \in \varpi_F^h \OF} \Bigl| \sum_{\underline{i}\in I_{\leq N}} a_{\underline{i}} \frac{z^{\underline{i}}}{\underline{i}!}\Bigr|.
\end{align}
\begin{proof}
La démonstration se fait par récurrence sur $N$, le résultat étant trivial si $N=0$. Supposons  $N>0$ et que l'inégalité \eqref{ineg} soit satisfaite pour tout $m < N$. Soient $P,Q \colon \OF \to E$ définies par:
\[
P(z) = \sum_{\underline{i}\in I_{\leq N}} a_{\underline{i}}  \frac{z^{\underline{i}}}{\underline{i}!} \quad \mbox{et} \quad  
Q(z) = \sum_{\underline{i}\in I_{\leq N-1}} a_{\underline{i}}  \frac{z^{\underline{i}}}{\underline{i}!}.
\]
Par la Proposition \ref{ultratec1} il existe une constante $C_1 \in \mathbb{R}_{\geq 1}$ et un $n_{0} \in \mathbb{Z}_{\geq 0}$ tels que pour tout $h \geq n_0$ et tout $\underline{m} \in I_{= N}$ on a:
\begin{align}\label{al}
|a_{\underline{m}}| q^{-Nh} \leq C_1  \sup_{z \in \varpi_F^h\OF} | P(z) |.
\end{align}
Il reste à majorer le terme $\sup_{\underline{i}\in I_{\leq N-1}} |a_{\underline{i}}|  q^{-h|\underline{i}|}$. Par hypothèse de récurrence il existe une constante $C_2 \in \mathbb{R}_{\geq 1}$ et un $n_0' \in \mathbb{Z}_{\geq 0}$ tels que pour tout $h \geq n_0'$ on a:
\begin{align}\label{diseq1111}
\sup_{\underline{i}\in I_{\leq N-1}} |a_{\underline{i}}| q^{-h|\underline{i}|} &\leq C_2 \sup_{z \in \varpi_F^h \OF} |Q(z)|. 
\end{align}
Remarquons que l'inégalité \eqref{al} et le fait que $|N!| \leqslant |\underline{i}!|$ pour tout $\underline{i} \in I_{\leq N}$ impliquent que pour tout $h\geq n_0$ et tout $z \in \varpi_F^h \OF$ on a:
\begin{align*}
\sup_{\underline{i} \in I_{=N}} \Big|a_{\underline{i}}  \frac{z^{\underline{i}}}{\underline{i}!} \Big| 
\leq  \sup_{\underline{i} \in I_{=N}} |a_{\underline{i}}| q^{-Nh}  |\underline{i}!| ^{-1} 
\leq C_1 \sup_{z \in \varpi_F^h\OF}| P(z) | |N!|^{-1}.
\end{align*}
et donc pour tout $h \geq n_0$ et tout $z \in \varpi_F^h \OF$ on obtient:
\begin{align}\label{diseq}
|Q(z)| = \Big| P(z) - \sum_{\underline{i}\in I_{=N}} a_{\underline{i}} \frac{z^{\underline{i}}}{\underline{i}!} \Big| 
 \leq C_1 \sup_{z \in \varpi_F^h\OF} | P(z)| |N!|^{-1}.
\end{align}
Posons $n_0'' = \sup\{n_0,n_0'  \}$. Alors par les inégalités \eqref{diseq1111} et \eqref{diseq} on a pour tout $h \geq n_0''$:  
\begin{align*}
\sup_{\underline{i}\in I_{\leq N-1}} |a_{\underline{i}}| q^{-h|\underline{i}|} & \leq C_1 C_2 \sup_{z \in \varpi_F^h \OF} |P(z)| |N!|^{-1} \\
& \leq C \sup_{z \in \varpi_F^h \OF} |P(z)|
\end{align*}
où l'on a posé $C = C_1 C_2 |N!|^{-1}$. Cela termine la preuve du lemme.
\end{proof}
\end{lemma}

\begin{prop}\label{grosso1}
Soit $f \in C^r(\OF,E)$. Alors
\begin{itemize}
\item[(i)] $D_{\underline{i}}f \in C^{r-|\underline{i}|}(\OF,E)$ pour tout $\underline{i} \in I_{\leq [r]}$. De plus il existe une constante $C \in \mathbb{R}_{\geq 0}$ telle que pour tout $\underline{i} \in I_{\leq [r]}$ on a:  
\[
\|D_{\underline{i}}f \|_{C^{r-|\underline{i}|}}  \leq C  \|f \|_{C^{r}}.
\] 
\item[(ii)] Soient $\underline{i}$ et $\underline{j}$ deux éléments de $I_{\leq [r]}$ tels que $|\underline{i}|+|\underline{j}|\leq [r]$. On a: 
\[
D_{\underline{j}}\big(D_{\underline{i}}f\big) =  D_{\underline{i}+\underline{j}}f.
\]
\end{itemize} 
\begin{proof}
Supposons $x \in \OF$ et $y,z \in \varpi_F^h \OF$. En développant $\sigma(y+z)^{i_{\sigma}}$ pour tout $\sigma \in S$ on a:
\begin{align*}
\sum_{\underline{i}\in I_{\leq [r]}}  D_{\underline{i}}f(x) \frac{(y+z)^{\underline{i}}}{\underline{i}!}  
&=  \sum_{\underline{i}\in I_{\leq [r]}}  D_{\underline{i}}f(x)  \sum_{\underline{k}\leqslant \underline{i}}  \frac{y^{\underline{k}}}{\underline{k}!} \, \frac{z^{\underline{i}-\underline{k}}}{(\underline{i}-\underline{k})!}   \\
&=  \sum_{\underline{i}\in I_{\leq [r]}} \sum_{\substack{
 \underline{k} \in \mathbb{Z}_{\geq 0}^{|S|} \\
|\underline{k}| \leqslant [r] - |\underline{i}|  
}} D_{\underline{k}+\underline{i}}f(x)  \frac{y^{\underline{k}}}{\underline{k}!} \, \frac{z^{\underline{i}}}{\underline{i}!}.
\end{align*}
On peut alors réécrire $\varepsilon_{f,[r]}(x,y+z)-\varepsilon_{f,[r]}(x+y,z)$ sous la forme
\begin{align}\label{maggio7}
-\sum_{\underline{i}\in I_{\leq [r]}} \frac{z^{\underline{i}}}{\underline{i}!} \Biggl(D_{\underline{i}}f(x+y) - \sum_{\substack{
 \underline{k} \in \mathbb{Z}_{\geq 0}^{|S|} \\
|\underline{k}| \leqslant [r] - |\underline{i}|  
}}  D_{\underline{k}+\underline{i}}f(x) \frac{y^{\underline{k}}}{\underline{k}!} \Biggr).
\end{align}
D'après les inégalités
\[
|\varepsilon_{f,[r]}(x,y+z)| q^{rh} \leq C_{f,r}(h)  \quad \mbox{et} \quad  |\varepsilon_{f,[r]}(x+y,z)| q^{rh} \leq C_{f,r}(h)
\]
on obtient la majoration
\begin{align}\label{maggio1}
| \varepsilon_{f,[r]}(x,y+z)-\varepsilon_{f,[r]}(x+y,z)| \leq C_{f,r}(h)  q^{-rh}.
\end{align}  
Notons $\varepsilon_{D_{\underline{i}}f,[r]}\colon \OF \times \OF \to E$, pour $\underline{i} \in I_{\leq [r]}$ et $h \in \mathbb{Z}_{\geq 0}$, la fonction définie par:
\begin{align}\label{eguaderiv}
 \varepsilon_{D_{\underline{i}}f,[r]}(x,y)  = D_{\underline{i}}f(x+y) - \sum_{\substack{
 \underline{k} \in \mathbb{Z}_{\geq 0}^{|S|} \\
|\underline{k}| \leqslant [r] - |\underline{i}|  
}}  D_{\underline{k}+\underline{i}}f(x) \frac{y^{\underline{k}}}{\underline{k}!}.
\end{align}
On veut montrer que pour tout $\underline{i} \in I_{\leq [r]}$ on a:
\[
\sup_{x \in \OF, y \in \varpi_F^h \OF}|\varepsilon_{D_{\underline{i}}f,[r]}(x,y)|  q^{h(r-|\underline{i}|)} \to 0 \quad \mbox{quand} \quad h \to +\infty.
\]
Fixons $x \in \OF$ et $y \in \varpi_F^h \OF$. Par le Lemme \ref{lemmatec} appliqué à $\sum_{\underline{i}\in I_{\leq [r]}} \varepsilon_{D_{\underline{i}}f,[r]}(x,y) \frac{z^{\underline{i}}}{\underline{i}!}$ il existe une constante $C \in \mathbb{R}_{\geq 0}$ et un $n_0$ tels que pour tout $\underline{i} \in I_{\leq [r]}$ et tout $l\geq n_0$ on a:
\begin{align}\label{maggio777}
|\varepsilon_{D_{\underline{i}}f,[r]}(x,y)| q^{-l|\underline{i}|} \leq C \sup_{z \in \varpi_F^l \OF}\Big|\sum_{\underline{i}\in I_{\leq [r]}} \varepsilon_{D_{\underline{i}}f,[r]}(x,y) \frac{z^{\underline{i}}}{\underline{i}!}  \Big|.
\end{align}
Les inégalités \eqref{maggio777} et \eqref{maggio1} impliquent alors que pour tout $h \geq n_0$ on a: 
\begin{align*}
\sup_{x\in \OF, y \in \varpi_F^h \OF} |\varepsilon_{D_{\underline{i}}f,[r]}(x,y)| q^{-h|\underline{i}|} \leq C  q^{-hr}  C_{f,r}(h)
\end{align*}
et donc pour tout $h \geq n_0$ on a:
\begin{align*}
\sup_{x\in \OF, y \in \varpi_F^h \OF}|\varepsilon_{D_{\underline{i}}f,[r]}(x,y)| q^{h(r-|\underline{i}|)} \leq C  C_{f,r}(h).
\end{align*}
On en déduit l'appartenance de $D_{\underline{i}}f$ à $C^{r-|\underline{i}|}(\OF,E)$ et la majoration 
\[
\|D_{\underline{i}}f \|_{C^{r-|i|}} \leq C   \|f \|_{C^{r}}.
\] 
Par ce qui précède la fonction $D_{\underline{j}}(D_{\underline{i}}f)$, pour $\underline{j}\in I_{\leq [r] - |\underline{i}|}$, est le terme de $\varepsilon_{D_{\underline{i}}f,[r]}(x,\cdot)$ qui est facteur de $\frac{y^{\underline{j}}}{\underline{j}!}$, c'est-à-dire 
\[
D_{\underline{j}}\big(D_{\underline{i}}f\big) =  D_{\underline{i}+\underline{j}}f,
\] 
d'où le résultat.
\end{proof}
\end{prop}

\begin{cor}\label{corollariocont}
L'application $D_{\underline{i}}$ de $C^{r}(\OF,E)$ dans $C^{r- |\underline{i}|}(\OF,E)$ définie par $f \mapsto D_{\underline{i}}f$ est bien définie et continue. 
\begin{proof}
Il s'agit d'une conséquence immédiate du (i) de la Proposition \ref{grosso1}.
\end{proof}  
\end{cor}

Rappelons qu'une \textit{$E$-algèbre de Banach} est une $E$-algèbre normée (c'est-à-dire vérifiant $\|ab\| \leq \|a\|\cdot \|b\|$ et $\|1\| \leq \|1\|$) telle que l'espace vectoriel normé sous-jacent soit un espace de Banach.

\begin{lemma}\label{prodottofunz}
Soit $r \in \mathbb{R}_{\geq 0}$. L'espace $C^r(\OF,E)$ est une $E$-algèbre de Banach. 
\begin{proof}
Soient $f,g \in C^r(\OF,E)$. Posons pour tout $x,y \in \OF$:
\[
\nu_{fg,[r]}(x,y) = f(x+y)g(x+y) - \Big(\sum_{\underline{i} \in I_{\leq [r]}} D_{\underline{i}}f(x) \frac{y^{\underline{i}}}{\underline{i}!}\Big) \Big(\sum_{\underline{j} \in I_{\leq [r]}} D_{\underline{j}}g(x) \frac{y^{\underline{j}}}{\underline{j}!}\Big).
\] 
Si on additionne et on soustrait le terme $f(x+y) \Big( \sum_{\underline{j} \in I_{\leq [r]}} D_{\underline{j}}g(x) \frac{y^{\underline{j}}}{\underline{j}!}  \Big)$, on obtient:
\begin{align}\label{dausare}
\nu_{fg,[r]}(x,y) = f(x+y) \varepsilon_{g,[r]}(x,y) + \Big(\sum_{\underline{j} \in I_{\leq [r]}} D_{\underline{j}}g(x) \frac{y^{\underline{j}}}{\underline{j}!}\Big) \varepsilon_{f,[r]}(x,y).
\end{align}
Comme $f$ est une fonction bornée sur $\OF$ et $g \in C^r(\OF,E)$ on a:
\begin{align}\label{algebbb1}
\sup_{x \in \OF, y \in \varpi_F^h \OF}|f(x+y) \varepsilon_{g,[r]}(x,y)| q^{rh} \to 0 \quad \mbox{quand} \quad h \to +\infty.
\end{align}
De manière analogue, comme les fonctions $D_{\underline{j}}g$, pour tout $\underline{j} \in I_{\leq [r]}$, sont bornées sur $\OF$ par définition et $g \in C^r(\OF,E)$ on a:
\begin{align}\label{algebbb2}
\sup_{x \in \OF, y \in \varpi_F^h \OF}\Big|\Big(\sum_{\underline{j} \in I_{\leq [r]}} D_{\underline{j}}g(x) \frac{y^{\underline{j}}}{\underline{j}!}\Big) \varepsilon_{f,[r]}(x,y) \Big| q^{rh} \to 0 \quad \mbox{quand} \quad h \to +\infty.
\end{align}
Et donc, par \eqref{algebbb1} et \eqref{algebbb2} on déduit 
\begin{align}\label{algeb2}
\sup_{x \in \OF, y \in \varpi_F^h \OF}|\nu_{fg,[r]}(x,y)| q^{rh} \to 0 \quad \mbox{quand} \quad h \to +\infty.
\end{align} 
Par ailleurs, on peut réécrire $\nu_{fg,[r]}(x,y)$ pour tout $x,y \in \OF$ sous la forme:
\begin{align*}
f(x+y)g(x+y) - \sum_{\underline{k} \in I_{\leq [r]}}\Big( \sum_{\underline{i}+ \underline{j} = \underline{k}}  \binom{\underline{k}}{\underline{i}} D_{\underline{i}}f(x) &D_{\underline{j}}g(x) \Big) \frac{y^{\underline{k}}}{\underline{k}!}+ \\ &+ \sum_{\underline{k} \in I_{> [r]}}\Big( \sum_{\underline{i}+ \underline{j} = \underline{k}}  \binom{\underline{k}}{\underline{i}} D_{\underline{i}}f(x) D_{\underline{j}}g(x) \Big) \frac{y^{\underline{k}}}{\underline{k}!}.
\end{align*}
Posons: 
\[
\forall x,y \in \OF,\quad \varepsilon_{fg,[r]}(x,y) = f(x+y)g(x+y) - \sum_{\underline{k} \in I_{\leq [r]}}\Big( \sum_{\underline{i}+ \underline{j} = \underline{k}}  \binom{\underline{k}}{\underline{i}} D_{\underline{i}}f(x) D_{\underline{j}}g(x) \Big) \frac{y^{\underline{k}}}{\underline{k}!} 
\]
Par \eqref{algeb2} et puisque l'on a:
\[
\sup_{x \in \OF, y \in \varpi_F^h \OF} \Big|\sum_{\underline{k} \in I_{> [r]}}\Big( \sum_{\underline{i}+ \underline{j} = \underline{k}}  \binom{\underline{k}}{\underline{i}} D_{\underline{i}}f(x) D_{\underline{j}}g(x) \Big) \frac{y^{\underline{k}}}{\underline{k}!}\Big| q^{rh} \to 0 \quad \mbox{quand} \quad h \to +\infty
\]
on déduit
\begin{align}\label{algeb}
\sup_{x \in \OF, y \in \varpi_F^h \OF}|\varepsilon_{fg,[r]}(x,y)| q^{rh} \to 0 \quad \mbox{quand} \quad h \to +\infty.
\end{align} 
La condition \eqref{algeb} implique que $fg$ est de classe $C^r$ et l'égalité suivante  
\[
\forall \underline{k} \in I_{\leq [r]}, \quad D_{\underline{k}}(fg) =  \sum_{\underline{i}+ \underline{j} = \underline{k}}  \binom{\underline{k}}{\underline{i}} D_{\underline{i}}f D_{\underline{j}}g.
\]
Il nous reste à montrer l'inégalité suivante:
\begin{align}\label{secondapa}
\forall f,g \in C^r(\OF,E), \quad  \|fg \|_{C^r} \leq \|f \|_{C^r} \|g \|_{C^r}.
\end{align}
\begin{itemize}
\item[(i)] On a:
\[
\sup_{x \in \OF} \Bigl|\frac{D_{\underline{k}}(fg)(x)}{\underline{k}!}\Bigr| \leq \sup_{\underline{i}+\underline{j} = \underline{k}} \Big( \sup_{x \in \OF} \Bigl|\frac{D_{\underline{i}}f(x)}{\underline{i}!}\Bigr| \sup_{x \in \OF} \Bigl|\frac{D_{\underline{j}}g(x)}{\underline{j}!}\Bigr| \Big) \leq \|f \|_{C^r} \|g \|_{C^r}.
\]
\item[(ii)] En utilisant l'égalité \eqref{dausare} on en déduit la minoration
\[
\sup_{x,y \in \OF} \frac{| \varepsilon_{fg,[r]}(x,y)|}{|y|^{r}} \leq  \|f \|_{C^r} \|g \|_{C^r}.
\]
\end{itemize}
En revenant à la définition de $\|\cdot \|_{C^r}$ on déduit de (i) et de (ii) l'inégalité \eqref{secondapa}.
\end{proof} 
\end{lemma}

\subsubsection{\textbf{Composition de fonctions}}
Soit $f\colon \OF \to E$ une fonction de classe $C^r$ et $h$ une fonction de $\OF$ dans $\OF$. Dans ce paragraphe nous donnons une condition suffisante sur $h$ pour que $f\circ h\colon \OF \to E$ soit de classe $C^r$. 

Commençon par la définition suivante.
\begin{defin}\label{definizione22}
Soit $r \in \mathbb{R}_{\geq 0}$. On dit que $h\colon \OF \to F$ est de classe $C^{r,id}$ sur $\OF$ s'il existe des fonctions bornées $h^{(i)} \colon \OF \to F$, pour $0\leq i \leq [r]$, telles que, si l'on définit $\varepsilon_{h,[r]}\colon \OF \times \OF \to F$ par: 
\begin{align*}
\varepsilon_{h,[r]}(x,y) = f(x+y) -
 \sum_{i = 0}^{[r]} h^{(i)}(x) \frac{y^{i}}{i!}
\end{align*}
et pour tout $k \in \mathbb{Z}_{\geq 0}$
\begin{align*}
C_{h,r}(k) = \sup_{x \in \OF, y \in \varpi_F^k \OF} |\varepsilon_{h,[r]}(x,y)| q^{rk}
\end{align*}
alors $C_{h,r}(k)$ tend vers $0$ quand $k$ tend vers $+\infty$. 
\end{defin}

Le même raisonnement donné dans la Remarque \ref{unifcont} nous dit qu'une fonction $h\colon \OF \to F$ de classe $C^{r,id}$ est continue. 

\begin{lemma}[Unicité]\label{unic22}
Soit $h$ une fonction de $\OF$ dans $F$ de classe $C^{r,id}$. Alors il existe une unique famille de fonctions 
\[
\big\{ h^{(i)} \colon \OF \to F, \ 0 \leq i \leq [r] \big\}
\]
qui vérifie la Définition \ref{definizione22}.
\begin{proof}
La preuve étant similaire à celle du Lemme \ref{unic}, nous y renvoyons le lecteur. 
\end{proof}
\end{lemma}

Notons $C^{r,id}(\OF,F)$ l'ensemble des fonctions de $\OF$ dans $F$ qui sont de classe $C^{r,id}$. On munit $C^{r,id}(\OF,F)$ de la norme $\|\cdot \|_{C^{r,id}}$ définie par:
\[
\|h \|_{C^{r,id}} = \sup \Big( \sup_{0 \leq i \leq [r]} \sup_{x \in \OF} \Bigl|\frac{h^{(i)}(x)}{i!}\Bigr|, \sup_{x,y \in \OF} \frac{| \varepsilon_{h,[r]}(x,y)|}{|y|^{r}}   \Big)
\]
ce qui en fait un espace de Banach sur $F$.

\begin{prop}[Composition de fonctions] \label{composizione}
Soit $r \in \mathbb{R}_{\geq 0}$. Si $h\colon \OF \to \OF$ est une fonction de classe $C^{r,id}$ alors
\begin{itemize}
\item[(i)] $\forall f \in C^r(\OF,E)$,  $f \circ h \in C^r(\OF,E)$;
\item[(ii)] l'application de $C^r(\OF,E)$ dans $C^r(\OF,E)$ définie par $f \mapsto f \circ h$ est continue.
\end{itemize}
  
\begin{proof}
Le cas $[r] = 0$ est facile et est laissé au lecteur. Supposons donc $[r]\geq 1$. Posons pour tout $x,y \in \OF$:
\begin{align}\label{formulacomp}
\varepsilon_{f \circ h,[r]}(x,y) = f(h(x+y)) - \sum_{\underline{j} \in I_{\leq [r]}} D_{\underline{j}}f(h(x)) \frac{\Big(\sum_{i = 1}^{[r]} h^{(i)}(x) \frac{y^{i}}{i!}\Big)^{\underline{j}}}{\underline{j}!}. 
\end{align} 
Si on additionne et on soustrait le terme $f\Big(h(x)+ \sum_{i=1}^{[r]} h^{(i)}(x) \frac{y^{i}}{i!}  \Big)$, on obtient:
\begin{align*}
f(h(x+y)) - f\Big(h(x)+ \sum_{i=1}^{[r]} h^{(i)}(x) \frac{y^{i}}{i!}  \Big) + f\Big(h(x) &+ \sum_{i=1}^{[r]} h^{(i)}(x) \frac{y^{i}}{i!}  \Big) - \\
&- \sum_{\underline{j} \in I_{\leq [r]}} D_{\underline{j}}f(h(x)) \frac{\Big(\sum_{i=1}^{[r]} h^{(i)}(x) \frac{y^{i}}{i!}\Big)^{\underline{j}}}{\underline{j}!}.
\end{align*}
On peut réécrire $f(h(x+y)) - f\Big(h(x)+ \sum_{i= 1}^{[r]} h^{(i)}(x) \frac{y^{i}}{i!}  \Big)$ sous la forme
\[
f(h(x+y))-f(h(x+y)- \varepsilon_{h,[r]}(x,y)).
\]
Comme $f$ est une fonction de classe $C^r$ sur $\OF$ et donc par le Lemme \ref{cl} de classe $C^1$, et $h \in C^{r,id}(\OF,F)$,  on a pour $y$ suffisamment petit: 
\begin{align}\label{lipi}
\bigl|f(h(x+y))-f(h(x+y)- \varepsilon_{h,[r]}(x,y))| \leq (\sup_{\sigma \in S}\sup_{x \in \OF}|D_{e_{\sigma}}f(x)| ) |\varepsilon_{h,[r]}(x,y)|.
\end{align}
L'inégalité \eqref{lipi} et le fait que $h$ est une fonction de classe $C^{r,id}$ sur $\OF$ impliquent 
\[
\sup_{x \in \OF, y \in \varpi_F^l \OF}|f(h(x+y))-f(h(x+y)- \varepsilon_{h,[r]}(x,y))|  q^{rl} \to 0 \quad \mbox{quand} \quad l \to +\infty. 
\]
Par ailleurs, comme $f$ est une fonction de classe $C^r$ sur $\OF$ alors
\begin{align*}
&\sup_{x \in \OF, y \in \varpi_F^l \OF} \Biggl|   f\Big(h(x)+\sum_{i = 1}^{[r]} h^{(i)}(x) \frac{y^{i}}{i!}  \Big) - \sum_{\underline{j} \in I_{\leq [r]}} D_{\underline{j}}f(h(x)) \frac{\Big(\sum_{i = 1}^{[r]} h^{(i)}(x) \frac{y^{i}}{i!}\Big)^{\underline{j}}}{\underline{j}!}\Biggr| q^{rl}   
\end{align*}
tend vers $0$ quand $l$ tend vers $+\infty$. Avec ce qui précède on en déduit:
\begin{align}\label{algeb1}
\sup_{x \in \OF, y \in \varpi_F^l \OF} |\varepsilon_{f \circ h,[r]}(x,y)|  q^{rl} \to 0 \quad \mbox{quand} \quad l \to +\infty.
\end{align} 
La condition \eqref{algeb1} implique que $f\circ h$ est une fonction de classe $C^r$ sur $\OF$ et que $D_{\underline{i}}(f\circ h)$ est le terme de $\varepsilon_{f \circ h,[r]}(\cdot,y)$ qui est facteur de $\frac{y^{\underline{i}}}{\underline{i}!}$. 

Il nous reste à montrer le (ii). 
\begin{itemize}  
\item[$\bullet$] Par ce qui précède on a:
\[
\forall \underline{i} \in I_{\leq [r]}, \quad \sup_{x \in \OF} \Bigl|\frac{D_{\underline{i}}(f \circ h)(x)}{\underline{i}!}\Bigr| \leq \sup(1, \|h \|_{C^{r,id}}^r)  \|f\|_{C^r}.  
\]
\item[$\bullet$] Soit $m\in \mathbb{Z}_{\geq 0}$ tel que $\varepsilon_{h,[r]}(x,y)\in \OF$ pour tout $x\in \OF$, $y\in \varpi_F^m\OF$. Par l'inégalité \eqref{lipi} on a:
\[
\sup_{x \in \OF,y\in \varpi_F^m\OF} \bigl|f(h(x+y))-f(h(x+y)- \varepsilon_{h,[r]}(x,y))|  |y|^{-r} \leq \|h\|_{C^{r,id}}\|f\|_{C^r}. 
\]
De plus, comme $f$ est une fonction de classe $C^r$ on a:
\begin{align*}
&\sup_{x \in \OF, y \in \varpi_F^m\OF} \Biggl|   f\Big(h(x)+\sum_{i = 1}^{[r]} h^{(i)}(x) \frac{y^{i}}{i!}  \Big) - \sum_{\underline{j} \in I_{\leq [r]}} D_{\underline{j}}f(h(x)) \frac{\Big(\sum_{i = 1}^{[r]} h^{(i)}(x) \frac{y^{i}}{i!}\Big)^{\underline{j}}}{\underline{j}!}\Biggr| |y|^{-r}  \\
\leq & \sup_{x \in \OF, y\in \varpi_F^m\OF} | \varepsilon_{f,[r]}(x,y) | |y|^{-r} \\
\leq & \ \ \|f \|_{C^r}.  
\end{align*}
\item[$\bullet$] Supposons $y\notin \varpi_F^m\OF$. En utilisant \eqref{formulacomp} on déduit
\[
\sup_{x\in \OF, y\notin \varpi_F^m\OF} \frac{|\varepsilon_{f\circ h,[r]}(x,y)|}{|y|^r} \leq q^{mr} \sup (\|f\|_{C^r},\|f\|_{C^r}\|h\|_{C^{r,id}}). 
\]  

\end{itemize}
En revenant à la définition de $\|\cdot \|_{C^r}$ on déduit qu'il existe une constante $C\in \mathbb{R}_{>0}$ telle que $\|f\circ h \|_{C^r} \leq C \|f\|_{C^r}$ pour tout $f\in C^r(\OF,E)$. D'où le résultat.

\end{proof}
\end{prop}

\section{Lien avec les fonctions localement analytiques}
\subsection{Espaces de fonctions localement $\Q$-analytiques}\label{ana}

Soit $J \subseteq S$ et $d_{\sigma} \in \mathbb{Z}_{\geq 0}$ pour $\sigma \in S\backslash J$.  Nous rappelons d'abord la construction de l'espace $\mathcal{F}(\OF, J, (d_{\sigma})_{\sigma \in S \setminus J})$ et puis nous montrons que cet espace s'injecte de façon continue dans $C^r(\OF,E)$. 

Pour $a \in \OF$ et $n \in \mathbb{Z}_{\geq 0}$, on note $\mathcal{O}(a+\uni^n\OF, J, (d_{\sigma})_{\sigma \in S \setminus J})$ le $E$-espace vectoriel des fonctions $f\colon a+\uni^n\OF \to E$ telles que 
\begin{align*}
f(z) = \sum_{\substack{
\underline{m} = (m_{\sigma})_{\sigma \in S} \in \mathbb{Z}_{\geq 0}^{|S|} \\
m_{\sigma} \leq d_{\sigma} \ \mbox{si} \ \sigma \in S\setminus J
}} a_{\underline{m}}(a) (z-a)^{\underline{m}}
\end{align*}
avec $a_{\underline{m}}(a) \in E$ et $|a_{\underline{m}}(a)| q^{-n(|\underline{m}|)} \to 0$ quand $|\underline{m}| \to +\infty$. Muni de la topologie définie par la norme  
\[
\| f \|_{a, n} = \sup_{\underline{m}} \left( |a_{\underline{m}}(a)| q^{-n (|\underline{m}|)} \right) 
\]
c'est un espace de Banach sur $E$. 

Pour $h \in \mathbb{Z}_{\geq 0}$ on note $\mathcal{F}_h(\OF, J, (d_{\sigma})_{\sigma \in S \setminus J})$ le $E$-espace vectoriel des fonctions $f\colon \OF \to E$ telles  que
\[
\forall a \in \OF, \quad f|_{a+\uni^{h}\OF} \in \mathcal{O}(a+\uni^{h}\OF, J, (d_{\sigma})_{\sigma \in S \setminus J}).
\]
On munit $\mathcal{F}_h(\OF, J, (d_{\sigma})_{\sigma \in S \setminus J})$ de la norme définie par
\begin{align}\label{formulanorma}
\| f \|_{\mathcal{F}_h} = \sup_{a\, \mathrm{mod}\, \uni^h} \|f|_{a+\uni^{h}\OF} \|_{a,h}
\end{align}
qui en fait un espace de Banach sur $E$. On voit immédiatement que cette définition ne dépend pas du choix des réprésentants. De plus (\cite[p. 107]{sch}) les inclusions 
\[
\mathcal{F}_h(\OF, J, (d_{\sigma})_{\sigma \in S \setminus J}) \subseteq \mathcal{F}_{h+1}(\OF, J, (d_{\sigma})_{\sigma \in S \setminus J})
\] 
sont continues et compactes.

\begin{defin}\label{definloc}
On note $\mathcal{F}(\OF, J, (d_{\sigma})_{\sigma \in S \setminus J})$ le $E$-espace vectoriel des fonctions $f\colon \OF \to E$ telles qu'il existe $h \in \mathbb{Z}_{\geq 0}$ tel que
\[
f \in \mathcal{F}_h(\OF, J, (d_{\sigma})_{\sigma \in S \setminus J}). 
\]
\end{defin}

On munit  l'espace $\mathcal{F}(\OF, J, (d_{\sigma})_{\sigma \in S \setminus J})$ de la topologie de la limite inductive qui en fait un espace de type compact.

\begin{prop}\label{grossaz}
Pour $h \in \mathbb{Z}_{\geq 0}$ et $r \in \mathbb{R}_{\geq 0}$ on a $\mathcal{F}_h(\OF,S) \subset C^r(\OF,E)$ et, quel que soit $f$ dans $\mathcal{F}_h(\OF,S)$ on a: 
\[
\|f \|_{C^r} \leq \|f \|_{\mathcal{F}_h} \, q^{rh}.
\]

\begin{proof}
Soit $f \in \mathcal{F}_h(\OF,S)$. Fixons un système de représentants $A_h$ des classes de $\OF/\varpi_F^h \OF$. Par hypothèse on peut écrire $f$ sous la forme
\[
f(x)= \sum_{
\underline{m} \in \mathbb{Z}_{\geq 0}^{|S|}
} a_{\underline{m}}(a)   (x-a)^{\underline{m}}
\]
quels que soient $a \in A_h$ et $x \in a+\uni^h\OF$. Posons pour tout $\underline{i} \in I_{\leq [r]}$, tout $a \in A_h$ et tout $x \in a+\uni^h\OF$:
\begin{align}\label{deriv}
\frac{D_{\underline{i}}f(x)}{\underline{i}!} =
\sum_{\substack{
\underline{m}  \in \mathbb{Z}_{\geq 0}^{|S|} \\
\underline{m}\geqslant \underline{i}
}} a_{\underline{m}}(a) \binom{\underline{m}}{\underline{i}} (x-a)^{\underline{m}-\underline{i}}.   
\end{align}
Montrons d'abord que  l'on a l'inégalité suivante: 
\begin{align}\label{maggio}
\forall x \in \OF, \quad \Big|\frac{D_{\underline{i}}f(x)}{\underline{i}!}\Big| \leq \|f \|_{\mathcal{F}_h} \, q^{h|\underline{i}|}.
\end{align}
En effet par la formule \eqref{deriv} on a les inégalités suivantes pour tout $a \in A_h$ et tout $x \in a+\uni^h\OF$:
\begin{align*}
\Big|\frac{D_{\underline{i}}f(x)}{\underline{i}!}\Big| &\leq \sup_{\underline{m}} \Big|a_{\underline{m}}(a) \binom{\underline{m}}{\underline{i}}(x-a)^{\underline{m}-\underline{i}} \Big| \\
&\leq \sup_{\underline{m}} \left( |a_{\underline{m}}(a)| q^{-h |\underline{m}|}\right) q^{h |\underline{i}|} \\
&= \|f\|_{\mathcal{F}_h} \, q^{h |\underline{i}|}.
\end{align*}
Comme le membre de droite ne dépend pas de $a \in A_h$ on en déduit l'inégalité pour tout $x \in \OF$. Posons:
\begin{align}\label{maggiort}
\forall x,y \in \OF, \quad \varepsilon_{f,[r]}(x,y) = f(x+y) - \sum_{\underline{i} \in I_{\leq [r]}} D_{\underline{i}}f(x) \frac{y^{\underline{i}}}{\underline{i}!}.
\end{align}
Donnons d'abord une estimation de $|\varepsilon_{f,[r]}(x,y)|$ pour tout $x \in \OF$ et tout $y \in \varpi_F^h \OF$. En développant $f(x+y)$ pour $x \in a+\uni^h\OF$ on obtient:
\begin{align*}
f(x+y) &= \sum_{\underline{m}\in \mathbb{Z}_{\geq 0}^{|S|}} a_{\underline{m}}(a)  \big((x-a)+y\big)^{\underline{m}} \\
&= \sum_{\underline{m} \in \mathbb{Z}_{\geq 0}^{|S|}} \sum_{\underline{i}\leqslant \underline{m}} a_{\underline{m}}(a)\binom{\underline{m}}{\underline{i}}(x-a)^{\underline{m}-\underline{i}} y^{\underline{i}}.
\end{align*} 
La formule \eqref{deriv} et celle ci-dessus permettent de réécrire $\varepsilon_{f,[r]}(x,y)$ sous la forme
\[
\varepsilon_{f,[r]}(x,y) = \sum_{\underline{i} \in I_{>[r]}} \frac{D_{\underline{i}}f(x)}{\underline{i}!}  y^{\underline{i}}.
\]
Par l'inégalité \eqref{maggio} on a:
\begin{align*}
|\varepsilon_{f,[r]}(x,y)| = \sup_{\underline{i} \in I_{>[r]}} \Big| \frac{D_{\underline{i}}f(x)}{\underline{i}!} y^{\underline{i}}  \Big| 
\leq \|f \|_{\mathcal{F}_h} \, q^{(h- val_F(y)/f) |\underline{i}|}  
\leq \|f \|_{\mathcal{F}_h} \, q^{(h- val_F(y)/f) ([r]+1)}   
\end{align*} 
(rappelons que $q = p^f$ est la cardinalité du corps résiduel de $F$) d'où
\begin{align}\label{maggios}
|\varepsilon_{f,[r]}(x,y)| q^{r val_F(y)/f} \leq \|f \|_{\mathcal{F}_h} \, q^{(h- val_F(y)/f) ([r]+1)}q^{r val_F(y)/f}.
\end{align}
Comme $-(val_F(y)/f)([r]+1) + r val_F(y)/f$ tend vers $-\infty$ quand $\mathrm{val}_F(y) \to +\infty$, cela montre que  $f$ est de classe $C^r$. 

La majoration $\|f \|_{C^r} \leq \|f \|_{\mathcal{F}_h}\, q^{rh}$ se déduit facilement de la formule \eqref{maggiort} et des majorations \eqref{maggio} et \eqref{maggios}.  
\end{proof}
\end{prop}

\begin{defin}
Notons $\mathcal{F}(\OF,S)$ l'espace $\mathcal{F}(U,J,(d_{\sigma})_{\sigma \in S \setminus J})$ pour $U = \OF$ et $J=S$ (et donc $S\backslash J = \emptyset$) ou, autrement dit, le $E$-espace vectoriel des fonctions $f\colon \OF \to E$ localement $\Q$-analytiques.   
\end{defin}

\begin{cor}\label{inicont}
L'espace $\mathcal{F}(\OF,S)$ s'injecte de façon continue dans $C^r(\OF,E)$.
\begin{proof}
Comme l'espace $\mathcal{F}(\OF,S)$ est muni de la topologie de la limite inductive le résultat est une conséquence immédiate de la Proposition \ref{grossaz}.
\end{proof}
\end{cor}

\subsection{Décomposition en vaguelettes des fonctions de classe $C^r$}\label{bannn}
Conservons les notations du §\ref{ana}. Si $I$ est un ensemble, on note $c_0(I)$ l'espace des fonctions $f\colon I \to E$ telles que:
\[
\forall \varepsilon > 0, \quad |\{ i \in I, |f(i)| \geq \varepsilon \}| < +\infty
\]
et on le munit de la norme définie par:
\[
\|f \|_{\infty} = \sup_{i \in I} |f(i)|.
\]
Soit $V$ un espace de Banach sur $E$. Une famille $(e_i)_{i \in I}$ d'éléments de $V$ est \textit{une base de Banach}  de $V$ si l'application de $c_0(I)$ dans $V$ définie par $(a_i)_{i \in I} \mapsto \sum_{i \in I} a_i e_i$ est un isomorphisme d'espaces de Banach. Rappelons que tout $E$-espace de Banach $V$ possède des bases de Banach (\cite[Proposition 10.1]{sch}) et que, si $\|V \| \subseteq |E|$, il existe un ensemble $I$ tel que $V$ est isométrique à $(c_0(I), \|\cdot \|_{\infty})$ (\cite[Remarque 10.2]{sch}).

Soit $r \in \mathbb{R}_{\geq 0}$. Le but de ce paragraphe est d'exhiber une agréable base de Banach de $C^r(\OF,E)$, qui dépend de $r$ et  qui consiste d'une famille dénombrable de fonctions localement $\Q$-polynomiales. Si $F = \Q$ cette base coïncide avec celle construite par  Van der Put pour l'espace des fonctions continues sur $\Z$ et généralisée par Colmez pour $r$ quelconque (\cite{vander}, \cite[Théorème I.5.14]{colmez}). Signalons que pour l'espace des fonctions continues sur $\OF$ cette base a déjà été construite par De Shalit \cite[§2]{sha}.

Posons $A_0 = \{0\}$. Choisissons pour tout $h \in \mathbb{Z}_{>0}$ un système de représentants $A_h \subset \OF$ des classes de  $\OF/\varpi_F^h \OF$ de sorte que $A_h \supset A_{h-1}$ et notons:
\[
A = \coprod_{h \geq 1} A_h \backslash A_{h-1}.
\] 
On a de même $A = \bigcup_{h \geq 1} A_h$. Pour tout $a \in A$ on définit $l(a)$ comme le plus petit entier $n_0$ tel que $a \in A_{n_0}$. De plus on fixe une numérotation $a_1, a_2, \ldots$ des éléments de $A$ de sorte que si $a_m \in A_i$ et $a_n \in A_{i+1} \backslash A_i$ alors $m < n$. 

Si $U$ est un ouvert de $\OF$ alors $\mathbf{1}_{U}$ désigne sa fonction caractéristique.

Posons pour tout $N,h \in \mathbb{Z}_{\geq 0}$: 
\begin{align*}
\mathcal{F}^{N}(\OF,S) & = \sum_{
\underline{d} \in I_{\leq N}} \mathcal{F}(\OF,\emptyset, \underline{d}); \\
\mathcal{F}_h^{N}(\OF,S) &= \sum_{
\underline{d} \in I_{\leq N}
} \mathcal{F}_h(\OF,\emptyset, \underline{d}). 
\end{align*} 
L'espace $\mathcal{F}^{N}(\OF,S)$ (resp. $\mathcal{F}_h^{N}(\OF,S)$) est  un sous-$E$-espace vectoriel de $\mathcal{F}(\OF,S)$ (resp. $\mathcal{F}_h(\OF,S)$).

\begin{prop}\label{put}
Soit $h \in \mathbb{Z}_{\geq 0}$.  
\begin{itemize}
\item[(i)] Les fonctions
\[
\mathbf{1}_{a+\varpi_F^{l(a)} \OF}(z)  \Big(\frac{z-a}{\varpi_F^{l(a)}}\Big)^{\underline{i}}
\]
pour $a \in A_h$ et $\underline{i}\in I_{\leq [r]}$ forment une base de $\mathcal{F}_h^{[r]}(\OF,S)$;
\item[(ii)] Les fonctions
\[
\mathbf{1}_{a+\varpi_F^{l(a)} \OF}(z)  \Big(\frac{z-a}{\varpi^{l(a)}}\Big)^{\underline{i}}
\]
pour $a \in A$ et $\underline{i}\in I_{\leq [r]}$ forment une base de $\mathcal{F}^{[r]}(\OF,S)$.
\end{itemize}
\begin{proof}
Notons $f_{a,\underline{i}}$ et $g_{a,\underline{i}}$ pour $a \in A_h$ et $\underline{i} \in I_{\leq [r]}$ les fonctions définies pour $z \in \OF$ par:
\begin{align*}
f_{a,\underline{i}}(z) &= \mathbf{1}_{a+\varpi_F^{h} \OF}(z)  \Big(\frac{z-a}{\varpi_F^{h}}\Big)^{\underline{i}}\\
g_{a,\underline{i}}(z) &= \mathbf{1}_{a+\varpi_F^{l(a)} \OF}(z)  \Big(\frac{z-a}{\varpi_F^{l(a)}}\Big)^{\underline{i}}.
\end{align*}
Les fonctions $f_{a,\underline{i}}$ pour $a \in A_h$ et $\underline{i}\in I_{\leq [r]}$ forment une base de $\mathcal{F}_h^{[r]}(\OF,S)$. On fixe deux numérotations: 
\begin{align*}
\eta &\colon A_h  \longrightarrow \mathbb{N} \\
\iota &\colon  A_h \times I_{\leq [r]}    \longrightarrow \mathbb{N}
\end{align*}  
de sorte que $\iota(a,\underline{i}) < \iota(b,\underline{j})$ si une parmi les trois conditions suivantes est satisfaite: 
\begin{itemize}
\item[(i)] $|\underline{i}| < |\underline{j}|$;
\item[(ii)] $|\underline{i}| = |\underline{j}|$ et $\underline{i} < \underline{j}$ par rapport à l'ordre lexicographique;
\item[(iii)] $\underline{i} = \underline{j}$ et $\eta(a) < \eta(b)$.
\end{itemize}
En écrivant $\Big(\frac{z-a}{\varpi_F^{l(a)}}\Big)^{\underline{i}} = \Big(\varpi_F^{h-l(a)} \, \frac{z-a-b\varpi_F^{l(a)}}{\varpi_F^h} + b\Big)^{\underline{i}}$ pour $z \in \OF$ et en développant on obtient: 
\begin{align*}
g_{a,\underline{i}}(z) = \mathbf{1}_{a+\varpi_F^{l(a)} \OF}(z) \Big(\frac{z-a}{\varpi_F^{l(a)}}\Big)^{\underline{i}} &= \Big( \sum_{b \in A_{h-l(a)}} \mathbf{1}_{a+b\varpi_F^{l(a)}+\varpi_F^h \OF}(z) \Big) \Big(\frac{z-a}{\varpi_F^{l(a)}}\Big)^{\underline{i}} \\
&=   \sum_{b \in A_{h-l(a)}} \sum_{
 \underline{m} \leqslant \underline{i} } s_{\underline{m}} f_{a+b\varpi_F^{l(a)},\underline{m}}(z)
\end{align*}
où l'on a noté
\[
s_{\underline{m}} =  \binom{\underline{i}}{\underline{m}} \big(\varpi_F^{h-l(a)}\big)^{\underline{m}} b^{\underline{i}-\underline{m}}.
\]
La matrice permettant de passer des $g_{a,\underline{i}}$ aux $f_{a,\underline{i}}$ pour $a \in A_h$ et $\underline{i} \in I_{\leq [r]}$ est triangulaire par blocs, à coefficients entiers, chaque bloc diagonal étant une matrice triangulaire inférieure avec des éléments de $\OF\backslash \{0\}$ sur la diagonale. Cela entraîne le (i) et le (ii). 
\end{proof}
\end{prop}

Fixons un plongement $\rho \colon F \into E$. Si $a \in A$ et $\underline{i} \in I_{\leq [r]}$, on note $e_{a,\underline{i},r}$ l'élément de $\mathcal{F}^{[r]}(\OF,S)$ défini par:
\[
z\mapsto e_{a,\underline{i},r} (z) = \unif^{[l(a)r]} \mathbf{1}_{a+\varpi_F^{l(a)} \OF}(z) \Big(\frac{z-a}{\varpi_F^{l(a)}}\Big)^{\underline{i}}.
\]

La Proposition \ref{put} implique la remarque suivante.
\begin{rem} \label{coefficienti}
{\rm(i) Les $e_{a,\underline{i},r}$ pour $a \in A_h$ et $\underline{i} \in I_{\leq [r]}$ forment une base de $\mathcal{F}_h^{[r]}(\OF,S)$.

(ii) Les $e_{a,\underline{i},r}$ pour $a \in A$ et $\underline{i} \in I_{\leq [r]}$ forment une base de $\mathcal{F}^{[r]}(\OF,S)$.}
\end{rem}

\begin{lemma} \label{lemval}
Si $\underline{i} \in I_{\leq [r]}$ alors
\[
\|e_{0,\underline{i},r} \|_{\Cr} = 1 \quad \mbox{et} \quad \|e_{a,\underline{i},r} \|_{\Cr} \leq q^{-([l(a)r]-l(a)r+r- |\underline{i}|)} \leq q \ \mbox{si} \ a\in A\backslash \{0 \}.
\]
\begin{proof}
Supposons $a = 0$ et notons $f = e_{0,\underline{i},r}$. Alors par l'Exemple \ref{esempio} on a: 
\[
 \sup_{
 \underline{m} \in I_{\leq [r]}} \sup_{x \in \OF} \Bigl|\frac{D_{\underline{m}}f(x)}{\underline{i}!}\Bigr| = 1 \quad \mbox{et} \quad \sup_{x,y \in \OF} \frac{| \varepsilon_{f,[r]}(x,y)|}{|y| ^{r}} = 0,
\]
ce qui entraîne $\|f \|_{\Cr} = 1$.

Supposons maintenant $a \in A\backslash \{0 \}$ (et donc $l(a) \geq 1$). Soit $\underline{i}\in I_{\leq [r]}$. Notons $g = e_{a,\underline{i},r}$. On a pour tout $\underline{m} \in I_{\leq [r]}$ tel que $\underline{m} \leqslant \underline{i}$:
\[
\frac{D_{\underline{m}}g(x)}{\underline{m}!} =
\unif^{[l(a)r]} \big(\varpi_F^{-l(a)} \big)^{\underline{m}} \mathbf{1}_{a+\varpi_F^{l(a)}\OF}(x) \binom{\underline{i}}{\underline{m}} \Big(\frac{x-a}{\varpi_F^{l(a)}} \Big)^{(\underline{i}-\underline{m})}   
\]
et $\frac{D_{\underline{m}}g(x)}{\underline{m}!}=0$ sinon. On en déduit pour tout $\underline{m} \in I_{\leq [r]}$:
\begin{align*}
\Big|\frac{D_{\underline{m}}g(x)}{\underline{m}!}\Big| &\leq q^{-[l(a)r]} q^{l(a)|\underline{m}|} 
\leq q^{-[l(a)r]} q^{l(a)|\underline{i}|} 
\leq q^{-([l(a)r]-l(a)r+r- |\underline{i}|)}
\end{align*} 
Notons que l'on a:
\[
\varepsilon_{g,[r]}(x,y) = \unif^{[l(a)r]} \Big(\mathbf{1}_{a+\varpi_F^{l(a)} \OF}(x+y) - \mathbf{1}_{a+\varpi_F^{l(a)} \OF}(x)   \Big) \Big(\frac{x+y-a}{\varpi_F^{l(a)}}\Big)^{\underline{i}},
\]
d'où $\varepsilon_{g,[r]}(x,y) = 0$ si $y \in \varpi_F^{l(a)}\OF$ ou si ni $x$ ni $x+y$ n'appartiennent à $a+\varpi_F^{l(a)}\OF$. Pour conclure il reste à étudier deux cas:
\begin{itemize}
\item[(i)] si $y \in \varpi_F^{l(a)-1}\OF \backslash \varpi_F^{l(a)}\OF$ et $x+y \in a+\varpi_F^{l(a)}\OF$, on a:
\[
\frac{| \varepsilon_{g,[r]}(x,y)|}{|y|^{r}} \leq q^{-[l(a)r]}  q^{r(l(a)-1)} \leq q^{-([l(a)r]-l(a)r+r- |\underline{i}|)}.
\]
\item[(ii)] si $y \in \varpi_F^{l(a)-1}\OF \backslash \varpi_F^{l(a)}\OF$ et $x \in a+\varpi_F^{l(a)}\OF$, on a:
\[
\frac{| \varepsilon_{g,[r]}(x,y)|}{|y|^{r}} = q^{-[l(a)r]}  q^{|\underline{i}|}  q^{r(l(a)-1)} \leq q^{-([l(a)r]-l(a)r+r- |\underline{i}|)}.
\]
\end{itemize}
En revenant à la définition de $\|\cdot \|_{C^r}$ on en déduit le lemme.
\end{proof}
\end{lemma}

\begin{lemma}\label{lemmatec2}
Soit $h \in \mathbb{Z}_{\geq 0}$. Si $b_a \in E$ pour tout $a \in A_h$ alors
\[
\sup_{a \in A_h} |b_a| \leq \big\|\sum_{a \in A_h} b_a e_{a,\underline{0},r}\big\|_{C^r} \, q^r.
\]
\begin{proof}
La démonstration se fait par récurrence sur $h$, le cas $h=0$ étant trivial. D'après la construction de l'ensemble $A$, pour tout $a \in A_h \backslash A_{h-1}$ il existe $u_a \in \OF^{\times}$ tel que $a + u_a \varpi_F^{h-1} \in A_{h-1}$. Notons $f$ la fonction définie par: 
\[
f(z)= \sum_{a \in A_h} b_a e_{a,\underline{0},r}(z).
\]
Comme les $e_{a,\underline{0},r}$ sont localement constantes on vérifie immédiatement que $D_{\underline{i}}f = 0$ pour tout $\underline{i} \in I_{\leq [r]}\backslash \{\underline{0}\}$ et donc on a $\varepsilon_{f,[r]}(x,y) = f(x+y)-f(x)$. De plus, pour tout $a \in A_{h-1}$ la fonction $e_{a,\underline{0},r}$ est constante modulo $\varpi_F^{h-1}\OF$ et, par construction de l'ensemble $A$, on a:  
\[
\sum_{c \in A_h\setminus A_{h-1}} b_c \mathbf{1}_{c+\varpi_F^{l(c)}\OF}(a) = \sum_{c \in A_h\setminus A_{h-1}} b_c \mathbf{1}_{c+\varpi_F^{h}\OF}(a) = 0.
\] 
En effet, si $a \in A_{h-1}$ alors $a \notin c + \varpi_F^h \OF$ quel que soit $c \in A_h$. On obtient alors pour $a \in A_h\backslash A_{h-1}$:
\[
\varepsilon_{f,[r]}(a,u_a  \varpi_F^{h-1}) = f(a+u_a \varpi_F^{h-1}) - f(a) = -\unif^{[rh]} b_a
\]
ce qui implique
\begin{align}\label{diseee}
\sup_{a \in A_h \setminus A_{h-1}} |b_a|   \leq  \sup_{a \in A_h \setminus A_{h-1}} q^{[hr]} |\varepsilon_{f,[r]}(a,u_a \varpi_F^{h-1})|  \leq q^{[hr]}  q^{r(1-h)}  \|f \|_{C^r},
\end{align}
d'où 
\begin{align}\label{primain}
\sup_{a \in A_h \setminus A_{h-1}} |b_a| \leq \|f \|_{C^r} \, q^r.
\end{align} 
De plus, par le Lemme \ref{lemval} si $a \in A_h\backslash A_{h-1}$ on a:
\begin{align}\label{diseee1}
\|e_{a,\underline{0},r} \|_{\Cr} \leq q^{-([hr]-hr+r)}.
\end{align} 
Par l'égalité
\[
\sum_{a \in A_{h-1}} b_a e_{a,\underline{0},r} = f- \sum_{a \in A_h\setminus A_{h-1}} b_a e_{a,\underline{0},r}
\]
couplée aux inégalités \eqref{diseee} et \eqref{diseee1} on déduit: 
\begin{align}\label{secondain}
\big\|\sum_{a \in A_{h-1}} b_a e_{a,\underline{0},r}\big\|_{C^r} \leq \sup \Big( \big\|f \big\|_{C^r}, \sup_{a \in A_h \setminus A_{h-1}} |b_a|  \big\|e_{a,\underline{0},r}\big\|_{C^r}\Big) \leq \|f \|_{C^r}.
\end{align}
Par hypothèse de récurrence on a:
\begin{align}\label{terzain}
\sup_{a \in A_{h-1}} |b_a| \leq \big\|\sum_{a \in A_{h-1}} b_a e_{a,\underline{0},r}\big\|_{C^r} \, q^r.
\end{align}
Les relations \eqref{primain}, \eqref{secondain} et \eqref{terzain} nous permettent alors de conclure.
\end{proof}
\end{lemma}

Si $f \in \mathcal{F}^{[r]}(\OF,S)$ on note $b_{a,\underline{i},r}(f)$, pour $a \in A$ et $\underline{i} \in I_{\leq [r]}$, les coefficients de $f$ dans la base des $e_{a,\underline{i},r}$.  

\begin{lemma}\label{lemmatec3}
Il existe $C \in \mathbb{R}_{\geq 0}$ tel que pour tout $f \in \mathcal{F}^{[r]}(\OF,S)$ on a:
\[
\sup_{a \in A} \sup_{\underline{i} \in I_{\leq [r]}} |b_{a,\underline{i},r}(f)| \leq C \|f\|_{C^r}.
\]
\begin{proof}
La démonstration se fait par récurrence sur $[r]$, le cas $[r] = 0$ étant trivial. Soit $\tau \in S$. Par hypothèse de récurrence il existe $C_{\tau} \in \mathbb{R}_{\geq 0}$ tel que: 
\[
\sup_{a \in A} \sup_{\substack{
\underline{i} \in I_{\leq [r]} \\
i_{\tau} \geq 1
}}
|b_{a,\underline{i},r}(f)| \leq C_{\tau} \|D_{e_{\tau}}f\|_{C^{r-1}}.
\]
D'après le (i) de la Proposition \ref{grosso1} il existe une constante $M \in \mathbb{R}_{\geq 0}$ telle que $\|D_{e_{\tau}}f\|_{C^{r-1}} \leq M \|f \|_{C^r}$ d'où l'inégalité
\[
\sup_{a \in A} \sup_{\substack{
\underline{i} \in I_{\leq [r]} \\
i_{\tau} \geq 1
}}
|b_{a,\underline{i},r}(f)| \leq C_{\tau}  M \|f\|_{C^{r}},
\]
et comme ceci vaut pour n'importe quel $\tau \in S$ on déduit
\begin{align}\label{intec1}
\sup_{a \in A} \sup_{\substack{
\underline{i} \in I_{\leq [r]} \\
|\underline{i}| \geq 1
}}
|b_{a,\underline{i},r}(f)| \leq \sup_{\tau \in S} \{C_{\tau} \}  M \|f\|_{C^{r}}.
\end{align}
Pour terminer, il ne reste qu'à majorer les coefficients restants. Considérons pour cela la fonction $f_0$ définie par:
\[
f_0 = f - \sum_{a \in A} \sum_{\substack{
\underline{i} \in I_{\leq [r]} \\
|\underline{i}| \geq 1
}}b_{a,\underline{i},r}(f) e_{a,\underline{i},r}. 
\]
Par le Lemme \ref{lemmatec2}, l'inégalité \eqref{intec1} et le Lemme \ref{lemval} on a:
\begin{align*}
\sup_{a \in A} |b_{a,\underline{0},r}(f)| \leq \|f_0 \|_{C^r} \, q^r \leq C  \|f \|_{C^r}
\end{align*}
où l'on a posé
\[
C = \sup \big\{q^r, \sup_{\tau \in S} \{C_{\tau} \} \, M \,  q^{r+1}  \big\}.
\]
D'où le résultat.
\end{proof}
\end{lemma}

\begin{prop}\label{densità}
Soit $h\in \mathbb{Z}_{\geq 0}$. Si $f \in C^r(\OF,E)$ on note $f_h$ l'élément de $\mathcal{F}^{[r]}(\OF,S)$ défini par:
\[
f_h(z) = \sum_{a \in A_h} \mathbf{1}_{a + \varpi_F^h \OF}(z) \Big( \sum_{\underline{i}\in I_{\leq [r]}} D_{\underline{i}}f(a) \, \frac{(z-a)^{\underline{i}}}{\underline{i}!} \Big).
\]
Alors: 
\begin{itemize}
\item[(i)] Il existe une constante $C \in \mathbb{R}_{\geq 0}$ et un $n_0 \in \mathbb{Z}_{\geq 0}$ tels que pour tout $h \geq n_0$ on a: 
\[
\|f_{h+1}-f_h \|_{\mathcal{F}_{h+1}} \leq C_{f,r}(h) C  q^{-rh};
\]
\item[(ii)] $f_h$ tend vers $f$ dans $C^r(\OF,E)$ quand $h$ tend vers $+\infty$.
\end{itemize}
\begin{proof}
(i) D'après l'égalité
\[
f_{h+1}(z) = \sum_{a \in A_h} \sum_{b \in A_1} \mathbf{1}_{a + b \varpi_F^h+ \varpi_F^{h+1} \OF}(z) \Big( \sum_{
 \underline{i} \in I_{\leq [r]}} D_{\underline{i}}f(a+b \varpi_F^h) \, \frac{(z-a-b \varpi_F^h)^{\underline{i}}}{\underline{i}!} \Big)
\]
on déduit que $f_{h+1}(z)-f_h(z)$ peut se récrire sous la forme
\begin{align*}
\sum_{a \in A_h} \sum_{b \in A_1} \mathbf{1}_{a + b\varpi_F^h+ \varpi_F^{h+1} \OF}(z) \Big( \sum_{
 \underline{i} \in I_{\leq [r]}} D_{\underline{i}}f(a+b \varpi_F^h) \, \frac{(z-a-b \varpi_F^h)^{\underline{i}}}{\underline{i}!}  
-D_{\underline{i}}f(a) \, \frac{(z-a)^{\underline{i}}}{\underline{i}!} \Big).
\end{align*}   
En développant, pour $a \in A_h$ et $b \in A_1$, $\sigma(z-a)^{i_{\sigma}} = \sigma((z-a-b \varpi_F^h)+b \varpi_F^h)^{i_{\sigma}}$ pour tout $\sigma \in S$ on a:
\begin{align*}
&\sum_{
 \underline{i} \in I_{\leq [r]}} D_{\underline{i}}f(a+b \varpi_F^h) \, \frac{(z-a-b \varpi_F^h)^{\underline{i}}}{\underline{i}!} 
-D_{\underline{i}}f(a) \, \frac{(z-a)^{\underline{i}}}{\underline{i}!}  \\
=& \sum_{
 \underline{i} \in I_{\leq [r]}} D_{\underline{i}}f(a+b \varpi_F^h) \, \frac{(z-a-b \varpi_F^h)^{\underline{i}}}{\underline{i}!} 
-\sum_{\substack{
 \underline{k} \in I_{\leq [r]} \\
\underline{k} \leqslant  \underline{i}
}}  D_{\underline{i}}f(a) \, \frac{(z-a-\varpi_F^h b)^{\underline{k}}}{\underline{k}!} \, \frac{(\varpi_F^h b)^{\underline{i}-\underline{k}}}{(\underline{i}-\underline{k})!}   \\
=&\sum_{
 \underline{j} \in I_{\leq [r]}} D_{\underline{j}}f(a+b \varpi_F^h) \, \frac{(z-a-b \varpi_F^h)^{\underline{j}}}{\underline{j}!} 
-\sum_{\substack{
 \underline{s} \in I_{\leq [r]} \\
|\underline{s}| \leqslant [r]- |\underline{j}|
}}  D_{\underline{j}+\underline{s}}f(a) \, \frac{(z-a-\varpi_F^h b)^{\underline{j}}}{\underline{j}!} \, \frac{(\varpi_F^h b)^{\underline{s}}}{\underline{s}!}  \\
=& \sum_{
 \underline{j} \in I_{\leq [r]}} \frac{1}{\underline{j}!} \Big( D_{\underline{j}}f(a+b \varpi_F^h)  
-\sum_{\substack{
 \underline{s} \in I_{\leq [r]} \\
|\underline{s}| \leqslant [r]- |\underline{j}|
}}  D_{\underline{j}+\underline{s}}f(a) \, \frac{(\varpi_F^h b)^{\underline{s}}}{\underline{s}!} \Big) (z-a-b \varpi_F^h)^{\underline{j}}.
\end{align*}
Par la preuve de la Proposition \ref{grosso1} il existe une constante $C \in \mathbb{R}_{\geq 0}$ et un $n_0 \in \mathbb{Z}_{\geq 0}$ tels que pour tout $h \geq n_0$, tout $a \in A_h$ et tout $b \in A_1$ on a:
\[
\Big|D_{\underline{j}}f(a+b \varpi_F^h)  
-\sum_{\substack{
 \underline{s} \in I_{\leq [r]} \\
|\underline{s}| \leq [r]- |\underline{j}|
}}  D_{\underline{j}+\underline{s}}f(a) \, \frac{(\varpi_F^h  b)^{\underline{s}}}{\underline{s}!} \Big| \leq q^{h (|\underline{j}|-r)}   C_{f,r}(h)  C.
\]
On en déduit pour tout $h \geq n_0$ 
\[
\|f_{h+1}-f_h \|_{\mathcal{F}_{h+1}} \leq \sup_{\underline{j}\in I_{\leq [r]}} \Big(q^{-(h+1)|\underline{j}|} q^{h (|\underline{j}|-r)}   C_{f,r}(h)  C \Big) = q^{-rh}   C_{f,r}(h)  C,
\]
d'où le (i).

(ii) D'après la Proposition \ref{grossaz} et par le (i) on a pour tout $h \geq n_0$:
\[
\|f_{h+1}-f_h \|_{C^r} \leq \|f_{h+1}-f_h \|_{\mathcal{F}_{h+1}} q^{r(h+1)} \leq q^r  C_{f,r}(h)  C 
\] 
et comme $C_{f,r}(h)$ tend vers $0$ quand $h$ tend vers $+\infty$ on en déduit que $f_h$ a une limite dans $C^r(\OF,E)$. Comme par ailleurs $f_h(a)=f(a)$ si $a \in A_h$ par définition de $f_h$, cette limite coïncide avec $f$ sur $A$, et donc partout par continuité, ce qui termine la preuve du (ii).
\end{proof}
\end{prop}

Le théorème suivant est essentiellement un corollaire des résultats que l'on a prouvés ci-dessus:

\begin{theo}\label{basebanach}
La famille des $e_{a,\underline{i},r}$, pour $a \in A$ et $\underline{i} \in I_{\leq [r]}$, forme une base de Banach de $C^r(\OF,E)$. 
\begin{proof}
Considérons l'application 
\[
\theta \colon c_0(A\times I_{\leq [r]}) \to C^r(\OF,E), \quad (b_{a,\underline{i}}) \mapsto \sum_{a \in A} \sum_{\underline{i}\in I_{\leq [r]}}b_{a,\underline{i}} e_{a,\underline{i},r}.
\]
Par le Lemme \ref{lemval} on déduit que $\theta$ est bien une application continue de $c_0(A\times I_{\leq [r]})$ dans $C^r(\OF,E)$. Notons $\varphi$ l'application de $\mathcal{F}^{[r]}(\OF,S)$ dans $c_0(A\times I_{\leq [r]})$ qui à tout $f \in \mathcal{F}^{[r]}(\OF,S)$ associe les coefficients de $f$ dans la base des $e_{a,\underline{i},r}$ (Remarque \ref{coefficienti} (ii)). Le Lemme \ref{lemmatec3} implique que cette application est continue une fois que l'on munit $\mathcal{F}^{[r]}(\OF,E)$ de la topologie induite par celle de $C^r(\OF,E)$. Comme $\mathcal{F}^{[r]}(\OF,S)$ est dense dans $C^r(\OF,E)$ par la Proposition \ref{densità}, alors $\varphi$ s'étend de façon unique en une application continue, que l'on désignera du même nom, de $C^r(\OF,E)$ dans $c_0(A\times I_{\leq [r]})$. Comme $\theta \circ \varphi (f) = f$ pour tout $f \in \mathcal{F}^{[r]}(\OF,E)$ on en déduit que $\theta \circ \varphi = id$. De manière analogue on a $  \varphi \circ  \theta = id$, ce qui prouve le résultat.
\end{proof}
\end{theo}

\begin{rem}
{\rm Soit $l \in \mathbb{R}_{\geq 0}$ tel que $l < r$. Alors $C^r(\OF,E) \subsetneq C^l(\OF,E)$. En effet il suffit de considérer la fonction $f$ définie par: 
\[
f = \sum_{\substack{
 \underline{i} \in I_{\leq [l]} \\
 a \in A  
}}\unif^{[r l(a)]-[l\, l(a)]}e_{a,\underline{i},l}.
\]
Le Théorème \ref{basebanach} implique que $f$ est une fonction de classe $C^l$ mais pas de classe $C^r$. }
\end{rem}

\subsection{Construction de sous-espaces fermés}\label{chiusi}
Soit $r \in \mathbb{R}_{\geq 0}$, $J \subseteq S$ et $d_{\sigma} \in \mathbb{Z}_{\geq 0}$ pour $\sigma \in S\backslash J$. Nous allons définir un sous-espace fermé de $C^r(\OF,E)$ qui dépend de $J$ et de $(d_{\sigma})_{\sigma \in S\setminus J}$ et puis nous construisons une base de Banach de cet espace à partir de la base de Banach de $C^r(\OF,E)$ décrite au §\ref{bannn}. 

Posons:
\[
J' = J \coprod \{\sigma \in S\backslash J,\, d_{\sigma}+1 > r\}
\] 
et désignons par $e_{\sigma}$ le vecteur de $\mathbb{Z}_{\geq 0}^{|S|}$ ayant toutes ses composantes nulles sauf celle d'indice $\sigma$ qui est égal à $1$. Notons pour tout $f \in C^r(\OF,E)$:   
\[
\forall \sigma \in S,\ 0\leq i \leq [r], \quad \frac{\partial^i }{\partial z_{\sigma}^i}f = D_{i e_{\sigma}}f.
\] 

\begin{defin}\label{sottosp}
On note $C^r(\OF,J',(d_{\sigma})_{\sigma \in S\setminus J'})$ le sous-$E$-espace vectoriel des fonctions $f$ de classe $C^r$ sur $\OF$ telles que 
\[
\forall \sigma \in S\backslash J', \quad \frac{\partial^{d_{\sigma}+1} }{\partial z_{\sigma}^{d_{\sigma}+1}}f = 0. 
\] 
\end{defin} 
  
D'après le Corollaire \ref{corollariocont}  l'opérateur $D_{\underline{i}}$ est continu pour tout $\underline{i} \in I_{\leq [r]}$ ce qui implique que l'espace $C^r(\OF,J',(d_{\sigma})_{\sigma \in S\setminus J'})$ est bien un sous-espace fermé de $C^r(\OF,E)$.  On le munit de la topologie induite par celle de $C^r(\OF,E)$ qui en fait un espace de Banach sur $E$.

Si $L$ désigne $J$ ou $J'$ posons pour tout $N \in \mathbb{Z}_{\geq 0}$: 
\[
\mathcal{F}^{N}(\OF,L, (d_{\sigma})_{\sigma \in S \setminus L}) = \mathcal{F}^{N}(\OF, S) \cap \mathcal{F}(\OF,L, (d_{\sigma})_{\sigma \in S \setminus L} ).
\] 
\begin{rem}\label{uguspazi}
{\rm Par définition de l'espace $\mathcal{F}^{[r]}(\OF,S)$ on voit facilement que:
\[
\mathcal{F}^{[r]}(\OF,J', (d_{\sigma})_{\sigma \in S \setminus J'}) = \mathcal{F}^{[r]}(\OF, J, (d_{\sigma})_{\sigma \in S \setminus J}).
\]
}
\end{rem}

Posons:
\begin{align*}
& Y = \bigl\{\underline{i}  \in \mathbb{Z}_{\geq 0}^{|S|}, \ i_{\sigma} \leq d_{\sigma} \ \mbox{si} \ \sigma \in S\backslash J \bigr\} \\
& Y' = \bigl\{\underline{i}  \in \mathbb{Z}_{\geq 0}^{|S|},\ i_{\sigma} \leq d_{\sigma} \ \mbox{si} \ \sigma \in S\backslash J' \bigr\},
\end{align*}
et notons que $Y\cap I_{\leq [r]} = Y' \cap I_{\leq [r]}$.

\begin{prop}\label{basebanach2}
La famille des $e_{a,\underline{i},r}$, pour $a \in A$ et $\underline{i} \in Y'\cap I_{\leq [r]}$, est une base de Banach de $C^r(\OF,J',(d_{\sigma})_{\sigma \in S\setminus J'})$. 
\begin{proof}
Soit $f \in C^r(\OF,J',(d_{\sigma})_{\sigma \in S\setminus J'})$ et $\tau \in S\backslash J'$. Par le Théorème \ref{basebanach} il existe une unique suite $(b_{a,\underline{i}})_{a,\underline{i}}$ d'éléments de $E$ tendant vers $0$ tel que $f$ s'écrit sous la forme
\[
f = \sum_{a \in A} \sum_{\underline{i}\in I_{\leq [r]}}b_{a,\underline{i}} e_{a,\underline{i},r}.
\]
Par continuité de l'opérateur $\frac{\partial^{d_{\tau}+1} }{\partial z_{\tau}^{d_{\tau}+1}}$ (Corollaire \ref{corollariocont}) on obtient:
\[
\frac{\partial^{d_{\tau}+1} }{\partial z_{\tau}^{d_{\tau}+1}}f = \sum_{a \in A} \sum_{\underline{i}\in I_{\leq [r]}} \frac{\partial^{d_{\tau}+1} }{\partial z_{\tau}^{d_{\tau}+1}}b_{a,\underline{i}} e_{a,\underline{i},r}.
\]
Or, si $\underline{j}$ désigne l'élément de $\mathbb{Z}_{\geq 0}^{|S|}$ défini par:
\[
j_{\sigma} =
\left\{ \begin{array}{ll}
i_{\sigma}-d_{\tau}-1   & \mbox{si} \ \sigma = \tau \\
i_{\sigma} & \mbox{sinon}
\end{array} \right.
\]
un calcul simple montre que
\[
\frac{\partial^{d_{\tau}+1} }{\partial z_{\tau}^{d_{\tau}+1}} e_{a,\underline{i},r} =
\left\{ \begin{array}{ll}
i_{\tau}(i_{\tau}-1)\ldots (i_{\tau}-d_{\tau}) e_{a,\underline{j},r}   & \mbox{si} \ i_{\tau}\geq d_{\tau}+1 \\
0 & \mbox{sinon}.
\end{array} \right.
\]
On en déduit que $\frac{\partial^{d_{\sigma}+1} }{\partial z_{\sigma}^{d_{\sigma}+1}}f = 0$ pour tout $\sigma \in S\backslash J'$ si et seulement si $f$ s'écrit sous la forme
\[
f = \sum_{a \in A}\ \sum_{
 \underline{i}\in Y'\cap I_{\leq [r]}  
} b_{a,\underline{i}} e_{a,\underline{i},r},
\]
d'où le résultat.
\end{proof}
\end{prop}

\begin{cor}\label{densità2}
Si $N$ est un entier tel que $N\geq [r]$ alors l'espace $\mathcal{F}^{N}(\OF,J, (d_{\sigma})_{\sigma \in S \setminus J})$ est dense dans  $C^r(\OF,J',(d_{\sigma})_{\sigma \in S\setminus J'})$.
\begin{proof}
La démonstration découle directement de la Proposition \ref{basebanach2}.
\end{proof}
\end{cor}

\section{Duaux}
Conservons les notations du  §\ref{chiusi} et notons $\mathcal{F}^{N}(\OF,J, (d_{\sigma})_{\sigma \in S \setminus J})^{\vee}$, pour tout $N\in \mathbb{Z}_{\geq 0}$,  l'ensemble des formes linéaires sur $\mathcal{F}^{N}(\OF,J, (d_{\sigma})_{\sigma \in S \setminus J})$. Si $N$ est un entier tel que $N\geq [r]$ alors, d'après le Corollaire \ref{densità2}, l'inclusion 
\[
\mathcal{F}^{N}(\OF,J,(d_{\sigma})_{\sigma \in S\setminus J}) \subseteq C^r(\OF,J',(d_{\sigma})_{\sigma \in S\setminus J'})
\] 
induit une injection
\[
C^r(\OF,J,(d_{\sigma})_{\sigma \in S\setminus J})^{\vee} \into \mathcal{F}^{N}(\OF,J,(d_{\sigma})_{\sigma \in S\setminus J})^{\vee}. 
\]
Dans cette section nous donnons une condition nécessaire et suffisante pour qu'une forme linéaire  $\mu\colon \mathcal{F}^{N}(\OF,J,(d_{\sigma})_{\sigma \in S\setminus J})\to E$ s'étende en une forme linéaire continue sur l'espace de Banach $C^r(\OF,J',(d_{\sigma})_{\sigma \in S\setminus J'})$. Cela généralise un résultat d\^u à Amice-Vélu et Vishik (\cite{amivel}, \cite{vis}).


\subsection{Distributions d'ordre $r$}\label{duali}
\begin{defin}
On appelle distribution tempérée d'ordre $r$ sur $\OF$ une forme linéaire continue sur l'espace de Banach $C^r(\OF,J',(d_{\sigma})_{\sigma \in S\setminus J'})$.
\end{defin} 

Notons:
\[
\big(C^r(\OF,J',(d_{\sigma})_{\sigma \in S \setminus J'})^{\vee}, \|\cdot \|_{\mathcal{D}_r,J',(d_{\sigma})_{\sigma}}\big)
\] 
l'espace des distributions tempérées d'ordre $r$ sur $\OF$ muni de la topologie forte.
 
Soit $N\in \mathbb{Z}_{\geq 0}$. Si $\mu \in \mathcal{F}^N(\OF,J,(d_{\sigma})_{\sigma \in S \setminus J})^{\vee}$ et $f \in \mathcal{F}^N(\OF,J,(d_{\sigma})_{\sigma \in S \setminus J})$ on note $\int_{\OF}f(z)\mu(z)$ l'accouplement et on pose:
\[
\int_{a+\uni^n\OF} f(z)\mu(z) = \int_{\OF} \mathbf{1}_{a+\uni^n\OF}(z) f(z) \mu(z)
\] 
où, pour $a \in \OF$ et $n \in \mathbb{Z}_{\geq 0}$,  $\mathbf{1}_{a+\uni^n\OF}$ désigne la fonction caractéristique de $a+\varpi_F^n \OF$.

\begin{theo}\label{velu}  (i) Soit $\mu \in C^r(\OF,J', (d_{\sigma})_{\sigma \in S \setminus J'})^{\vee}$. Il existe une constante $C_{\mu} \in \mathbb{R}_{\geq 0}$ telle que pour tout $a \in \OF$, tout $n \in \mathbb{Z}_{\geq 0}$ et tout  $\underline{i} \in Y'$ on ait: 
\begin{align}\label{inte1}
\Big|\int_{a+\uni^n\OF} \Big(\frac{z-a}{\varpi_F^n} \Big)^{\underline{i}} \mu(z)\Big| \leq C_{\mu} \, q^{nr}. 
\end{align}

(ii) Soit $N$ un entier tel que $N\geq [r]$ et $\mu \in  \mathcal{F}^N(\OF,J, (d_{\sigma})_{\sigma \in S \setminus J})^{\vee}$. Supposons qu'il existe une constante $C_{\mu} \in \mathbb{R}_{\geq 0}$ telle que pour tout $a \in \OF$, tout $n \in \mathbb{Z}_{\geq 0}$ et tout $\underline{i} \in Y\cap I_{\leq N}$ on ait:    
\begin{align}\label{inte2}
\Big| \int_{a+\uni^n\OF}  \Big(\frac{z-a}{\varpi_F^n} \Big)^{\underline{i}} \mu(z)\Big| \leq C_{\mu} \, q^{nr}.
\end{align}
Alors $\mu$ se prolonge de manière unique en une distribution tempérée d'ordre $r$ sur $\OF$. 
\begin{proof}
(i) Soient $a \in \OF$, $n \in \mathbb{Z}_{\geq 0}$ et $\underline{i} \in Y'$.  Notons $f_{a,\underline{i},n}$ la fonction définie par:
\[
f_{a,\underline{i},n}(z) = \mathbf{1}_{a+\uni^n\OF}(z) \Big(\frac{z-a}{\varpi_F^n} \Big)^{\underline{i}}.
\]
C'est un élément de l'espace $\mathcal{F}_{n}(\OF,S)$. Rappelons que $\mathcal{F}_{n}(\OF,S)$ est un $E$-espace de Banach, la norme $\|\cdot \|_{\mathcal{F}_n}$ étant définie par la formule \eqref{formulanorma} (§\ref{ana}).  On vérifie facilement que $\|f_{a,\underline{i},n} \|_{\mathcal{F}_n} = 1$. En utilisant la Proposition \ref{grossaz} on obtient:
\begin{align*}\label{normacr}
| \mu(f_{a,\underline{i},n})| &\leq \|\mu \|_{\mathcal{D}_r,J',(d_{\sigma})_{\sigma}}  \|f_{a,\underline{i},n} \|_{C^r} \\
&\leq \|\mu \|_{\mathcal{D}_r,J',(d_{\sigma})_{\sigma}}  q^{rn} \|f_{a,\underline{i},n} \|_{\mathcal{F}_n} \\
&= \|\mu \|_{\mathcal{D}_r,J',(d_{\sigma})_{\sigma}}  q^{rn}
\end{align*}
d'où le résultat une fois que l'on a posé $C_{\mu} = \|\mu \|_{\mathcal{D}_r,J',(d_{\sigma})_{\sigma}}  q^{rn}$.

(ii) L'unicité d'une telle extension découle du Corollaire \ref{densità2}: l'espace  $\mathcal{F}^{N}(\OF,J, (d_{\sigma})_{\sigma \in S \setminus J})$ est dense dans $C^r(\OF,J', (d_{\sigma})_{\sigma \in S \setminus J'})$. Montrons l'existence. Si $a \in A$ et $\underline{i} \in I_{\leq [r]} \cap Y$ posons:
\[
b_{a,\underline{i}} = \mu (e_{a,\underline{i},r} ) = \mu \Big(\unif^{[rl(a)]} \mathbf{1}_{a+\uni^{l(a)}\OF}(z)  \Big(\frac{z-a}{\varpi_F^{l(a)}} \Big)^{\underline{i}}  \Big). 
\]  
Comme par hypothèse il existe une constante $C_{\mu} \in \mathbb{R}_{\geq 0}$ telle que pour tout $a \in \OF$, tout $n \in \mathbb{Z}_{\geq 0}$ et tout $\underline{i} \in Y\cap I_{\leq N}$ on a:     
\begin{align*}
\Big| \mu \Big(\mathbf{1}_{a+\uni^n\OF}(z) \Big(\frac{z-a}{\varpi_F^n} \Big)^{\underline{i}} \Big) \Big| \leq C_{\mu} \, q^{nr}, 
\end{align*}
on déduit $|b_{a,\underline{i}}| \leq C_{\mu} \, q$ pour tout  $a \in A$ et tout  $\underline{i} \in I_{\leq [r]} \cap Y$. Par la Proposition \ref{basebanach2} il existe un unique élément $\tilde{\mu}$ de $C^r(\OF,J', (d_{\sigma})_{\sigma \in S \setminus J'})^{\vee}$ tel que l'on a $\tilde{\mu}( e_{a,\underline{i},r} ) = b_{a,\underline{i}}$ pour tout $a \in A$ et tout $\underline{i} \in I_{\leq [r]} \cap Y$. Notons: 
\[
\lambda = \tilde{\mu}|_{\mathcal{F}^N(\OF,J, (d_{\sigma})_{\sigma \in S \setminus J})} - \mu.
\] 
Montrons que $\lambda$ est identiquement nulle sur $\mathcal{F}^N(\OF,J, (d_{\sigma})_{\sigma \in S \setminus J})$. Remarquons qu'elle est identiquement nulle sur $\mathcal{F}^{[r]}(\OF,J, (d_{\sigma})_{\sigma \in S \setminus J})$ par construction. De plus, le point (i) implique que $\tilde{\mu}|_{\mathcal{F}^N(\OF,J, (d_{\sigma})_{\sigma \in S \setminus J})}$ satisfait \eqref{inte2}. Cela implique que $\lambda$ satisfait aussi \eqref{inte2}. Soit $a \in \OF$, $\underline{i} \in Y \cap I_{N}$ et $m \in \mathbb{Z}_{\geq 0}$. On peut réécrire la fonction $\mathbf{1}_{a+\uni^n\OF}(z) z^{\underline{i}}$ sous la forme:
\[
\sum_{b \in A_m} \sum_{\underline{s}\leqslant \underline{i}} \mathbf{1}_{(a+b\varpi_F^n)+  \uni^{n+m}\OF}(z)  \binom{\underline{i}}{\underline{s}} (a+b \varpi_F^n)^{\underline{i}-\underline{s}} (z-(a+b \varpi_F^n))^{\underline{s}}.
\]
Or, si $|\underline{s}| \leq [r]$ on a:
\[
\lambda \Big(\mathbf{1}_{(a+b\varpi_F^n)+\uni^{n+m}\OF}(z)(z-(a+b \varpi_F^n))^{\underline{s}}\Big) = 0. 
\]
Supposons donc $|\underline{s}| > [r]$. Comme $\lambda$ satisfait \eqref{inte2} déduit:
\[
\Big|\lambda \Big(\mathbf{1}_{(a+b\varpi_F^n)+\uni^{n+m}\OF}(z) (z-(a+b \varpi_F^n))^{\underline{s}}\Big) \Big| \leq C q^{(n+m)(r-|\underline{s}|)}
\]
d'où
\[
\Big|\lambda \Big(\mathbf{1}_{a+\uni^n\OF}(z) z^{\underline{i}}\Big) \Big| \leq C q^{(n+m)(r-([r]+1))}.
\]
On voit facilement que cette quantité tend vers $0$ quand $m$ tend vers $+\infty$, et cela implique
\[
\lambda \Big(\mathbf{1}_{a+\uni^n\OF}(z) z^{\underline{i}} 
 \Big) = 0.
\] 
On en déduit l'égalité $\tilde{\mu}|_{\mathcal{F}^N(\OF,J, (d_{\sigma})_{\sigma \in S \setminus J})} = \mu$, ce qui permet de conclure. 
\end{proof}
\end{theo}

Par ce qui précède on peut munir l'espace $C^r(\OF,J', (d_{\sigma})_{\sigma \in S \setminus J'})^{\vee}$ d'une norme équivalente à la norme $\|\cdot \|_{\mathcal{D}_r,J',(d_{\sigma})_{\sigma}}$, mais plus commode. 

\begin{cor}\label{normaequiv}
(i) Si l'on définit $\|\mu \|_{r,Y}$, pour $\mu \in C^r(\OF,J', (d_{\sigma})_{\sigma \in S \setminus J'})^{\vee}$ par la formule
\[
\|\mu \|_{r,Y} = \sup_{a \in \OF, n \in \mathbb{Z}_{\geq 0}} \sup_{\underline{i}\in Y} \Big( \Big| \int_{a+\uni^n\OF} \Big(\frac{z-a}{\varpi_F^n} \Big)^{\underline{i}} \mu(z)\Big| q^{-nr} \Big)
\] 
alors $\|\cdot \|_{r,Y}$ est une norme sur $C^r(\OF,J', (d_{\sigma})_{\sigma \in S \setminus J'})^{\vee}$ équivalente à  $\|\cdot \|_{\mathcal{D}_r,J',(d_{\sigma})_{\sigma}}$.

(ii) Si $N\geq [r]$ et si l'on définit $\|\mu \|_{r,Y}$, pour $\mu \in C^r(\OF,J', (d_{\sigma})_{\sigma \in S \setminus J'})^{\vee}$ par la formule
\[
\|\mu \|_{r,N} = \sup_{a \in \OF, n \in \mathbb{Z}_{\geq 0}} \sup_{\underline{i}\in Y\cap I_{\leq N}} \Big( \Big| \int_{a+\uni^n\OF} \Big(\frac{z-a}{\varpi_F^n} \Big)^{\underline{i}} \mu(z)\Big| q^{-nr} \Big)
\] 
alors $\|\cdot \|_{r,N}$ est une norme sur $C^r(\OF,J', (d_{\sigma})_{\sigma \in S \setminus J'})^{\vee}$ équivalente à  $\|\cdot \|_{\mathcal{D}_r,J',(d_{\sigma})_{\sigma}}$.
\begin{proof}
Les preuves de (i) et (ii) étant similaires, on se contentera de prouver la première assertion.  Par le Théorème \ref{velu}, on déduit facilement que $\|\cdot \|_{r,Y}$ est une norme sur l'espace $C^r(\OF,J', (d_{\sigma})_{\sigma \in S \setminus J'})^{\vee}$. De plus, l'application identité:
\[
id \colon (C^r\big(\OF,J', (d_{\sigma})_{\sigma \in S \setminus J'})^{\vee},\|\cdot \|_{\mathcal{D}_r,J',(d_{\sigma})_{\sigma}}\big) \to \big(C^r(\OF,J', (d_{\sigma})_{\sigma \in S \setminus J'})^{\vee}, \|\cdot \|_{r,Y}\big)
\]
est continue par l'inégalité de la preuve de (i) du Théorème \ref{velu}. Et donc, d'après le théorème de l'image ouverte (\cite[Proposition 8.6]{sch}), c'est un isomorphisme de $E$-espaces de Banach ce qui implique que la norme $\|\cdot \|_{\mathcal{D}_r,J',(d_{\sigma})_{\sigma}}$ est équivalente à la norme $\|\cdot \|_{r,Y}$.
 
\end{proof}
\end{cor}

\section{Une autre notion de fonction de classe $C^r$ sur $\OF$}

Soit $r \in \mathbb{R}_{\geq 0}$. Conservons les notations des sections §\ref{classe}, §\ref{duali} et §\ref{ana}, et notons $d = [F:\Q]$. En utilisant le fait que $\OF$ est un $\Z$-module libre de rang $d$ on est amené à considérer une autre notion, tout à fait naturelle, de fonction de classe $C^r$ sur $\OF$. Le but de cette section est de montrer que cette nouvelle notion n'est pas équivalente à celle donnée au §\ref{classe} dès que $r > 0$. 

Fixons un  $d$-uplet $\vec{r}= (r_i)_{1 \leq i \leq d}$ de nombres réels positifs ou nuls tels que $\sum r_i = r$ et une base $(e_i)_{1 \leq i \leq d}$ de $\OF$ sur $\Z$.  Notons $\theta$  l'isomorphisme de $\Z$-modules défini par:
\[
\theta\colon \Z^d \stackrel{\sim}{\longrightarrow} \OF, \quad (a_1, \ldots, a_d) \mapsto \sum_{i=1}^d a_i e_i.
\]
Si $z \in \bigotimes_{i=1}^d C^{r_i}(\Z,E)$, on définit $\|z \|$ comme l'infimum des $\sup_{j\in J} \|v_{j_1}\|_{C^{r_1}}\cdot \ldots \cdot \|v_{j_d}\|_{C^{r_d}}$ pour toutes les écritures possibles de $z$ sous la forme $\sum_{j\in J} v_{j_1}\otimes \ldots \otimes v_{j_d}$. Ceci munit $\bigotimes_{i=1}^d C^{r_i}(\Z,E)$ d'une semi-norme et on note $\widehat{\bigotimes}_{i=1}^d C^{r_i}(\Z,E)$ le séparé complété de l'espace $\bigotimes_{i=1}^d C^{r_i}(\Z,E)$ pour cette semi-norme. 

Posons:
\[
C^{\vec{r}}(\OF,E) = \big\{f\colon \OF \to E, \ \textstyle{f\circ \theta \in  \widehat{\bigotimes}_{i=1}^d C^{r_i}(\Z,E)}  \big\},
\]
et munissons $C^{\vec{r}}(\OF,E)$ de la topologie déduite de celle définie sur $\widehat{\bigotimes}_{i=1}^d C^{r_i}(\Z,E)$. Notons $C^{\vec{r}}(\OF,E)^{\vee}$ le dual continu de l'espace de Banach $C^{\vec{r}}(\OF,E)$ muni de la topologie forte.

\begin{rem}
{\rm L'application de $\bigotimes_{i=1}^d \C^{0}(\Z,E)$ dans $C^0(\Z^d,E)$ définie par:
\[
\psi_1 \otimes \ldots \otimes\psi_d \longmapsto [(z_1,\ldots,z_d)\mapsto \psi_1(z_1)\cdot \ldots \cdot \psi_d(z_d)]
\]
s'étende en un isomorphisme de $E$-espaces de Banach:
\begin{align} \label{contin}
\textstyle{  \widehat{\bigotimes}_{i=1}^d C^{0}(\Z,E)} \stackrel{\sim}{\rightarrow} C^0(\Z^d,E)
\end{align}
(voir \cite[§17]{sch}). Or, comme $C^0(\Z^d,E)$ est isomorphe topologiquement à $C^0(\OF,E)$, on déduit de \eqref{contin} que $C^0(\OF,E)$ est isomorphe topologiquement à $C^{\vec{0}}(\OF,E)$.  }
\end{rem}

Supposons $r > 0$. Posons:
\[
\mathrm{Pol}(\OF,E) = \bigoplus_{n \in \mathbb{Z}_{\geq 0}} \mathcal{F}^n (\OF,S)
\]
où 
\[
\forall n \in \mathbb{Z}_{\geq 0}, \quad \mathcal{F}^{n}(\OF,S)  = \sum_{
\underline{d} \in I_{\leq n}} \mathcal{F}(\OF,\emptyset, (d_{\sigma})_{\sigma \in S}),
\]  
et notons que $\theta$ induit un isomorphisme de $E$-espaces vectoriels
\begin{align}\label{isopoli}
\mathrm{Pol}(\OF,E) \stackrel{\sim}{\rightarrow} \textstyle{ \bigotimes_{i = 1}^d \mathrm{Pol}(\Z,E)}, 
\end{align}
où $\mathrm{Pol}(\Z,E)$ désigne le $E$-espace vectoriel des fonction $f\colon \Z \to E$ localement polynomiales. 
 
Posons $B_0 = \{\vec{0} \}$. Choisissons pour tout $h \in \mathbb{Z}_{> 0}$ un système de représentants $B_h \subset \Z^d$ des classes de $\Z^d/p^h \Z^d$ de sorte que $B_h \supset B_{h-1}$. 
 
\begin{prop}\label{noniso}
Les espaces de Banach $C^r(\OF,E)^{\vee}$ et $C^{\vec{r}}(\OF,E)^{\vee}$ ne sont pas isomorphes.
\begin{proof} 
Soit $\mu \in \mathrm{Pol}(\OF,E)^{\vee}$. D'après le Théorème \ref{velu} on sait que $\mu$ s'étende à $C^r(\OF,E)$ si et seulement s'il existe une constante $C_{\mu} \in \mathbb{R}_{\geq 0}$ telle que pour tout $a \in \OF$, tout $n \in \mathbb{Z}_{\geq 0}$ et tout  $\underline{i} \in \mathbb{Z}_{\geq 0}^d$ on a: 
\begin{align*}
\Big|  \mu\Big(\mathbf{1}_{a+\uni^n\OF}(z) \Big(\frac{z-a}{\varpi_F^n} \Big)^{\underline{i}}\Big)\Big| \leq C_{\mu} \, q^{nr} 
\end{align*} 
ou, de manière équivalente en utilisant l'isomorphisme \eqref{isopoli}, si et seulement s'il existe une constante $C_{\mu} \in \mathbb{R}_{\geq 0}$ telle que pour tout $(a_i)_i \in \Z^d$, tout $n \in \mathbb{Z}_{\geq 0}$ et tout  $(j_i)_i \in \mathbb{Z}_{\geq 0}^d$ on a:    
\begin{align}\label{primafac}
\Big|\mu\Big(\mathbf{1}_{a+p^n \Z^d}(z) \prod_{i=1}^d \Big(\frac{z_i-a_i}{p^n} \Big)^{j_i} \Big)   \Big| \leq C_{\mu} q^{enr},
\end{align}
où $e$ désigne l'indice de ramification de $F$.

D'autre part, une application immédiate du critère d'Amice-Vélu et Vishik (\cite{amivel}, \cite{vis}) montre que la distribution $\mu$ s'étende à $C^{\vec{r}}(\OF,E)$ si et seulement s'il existe une constante $C_{\mu} \in \mathbb{R}_{\geq 0}$ telle que pour tout $(a_i)_i \in \Z^d$, tout $(n_i)_i \in \mathbb{Z}_{\geq 0}^d$ et tout $(j_i)_i \in \mathbb{Z}_{\geq 0}^d$ on a:  
\begin{align}\label{secondafac}
\Big|\mu\Big(\mathbf{1}_{\prod_i (a_i +p^{n_i} \Z^d)}(z) \prod_{i=1}^d \Big(\frac{z_i-a_i}{p^{n_i}} \Big)^{j_i} \Big)    \Big| \leq C_{\mu} q^{\sum_{i} e n_i r_i}.
\end{align}

Nous construisons maintenant une distribution $\mu \in \mathrm{Pol}(\OF,E)^{\vee}$ qui vérifie  \eqref{primafac} mais qui ne vérifie pas \eqref{secondafac}. Rappelons qu'un élément de $\mathrm{Pol}(\OF,E)^{\vee}$ est équivalent à la donnée des valeurs:
\[
\mu\Big(\mathbf{1}_{a+p^n \Z^d}(z) \prod_{i=1}^d z_i^{j_i} \Big)
\] 
pour $a \in \Z^d$, $(j_i)_i \in \mathbb{Z}_{\geq 0}^d$ et $n \in \mathbb{Z}_{\geq 0}$ avec les relations de compatibilité évidentes:
\begin{align}\label{reladicomp}
\mu\Big(\mathbf{1}_{a+p^n \Z^d}(z) \prod_{i=1}^d z_i^{j_i} \Big) = \sum_{k \in I_{\leq p-1}} \mu\Big(\mathbf{1}_{a+kp^n+p^{n+1} \Z^d}(z) \prod_{i=1}^d z_i ^{j_i} \Big).
\end{align}
Par hypothèse il existe $k$, $1\leq k \leq d$ tel que $r_k < r$. Pour $n\in \mathbb{Z}_{>0}$ notons $X_n = \prod_i X_{i,n}$ l'ouvert de $\Z^d$ défini par:
\[
X_{j,n} = \left\{ \begin{array}{ll}
\Z \quad &\mbox{si} \ j=k   \\
p^n\Z \quad &\mbox{si} \ j\neq k.
\end{array} \right.
\] 
Fixons $\alpha \in B_1\backslash X_1$. Pour tout $(j_i)_i \in \mathbb{Z}^d$ posons:
\[
\mu\Big(\mathbf{1}_{\Z^d}(z) \prod_{i=1}^d z_i^{j_i}\Big) = 1,
\]
et définissons $\mu$ sur les ouverts de la forme $b+p^n \Z^d$ pour $n \geq 1$ et $b \in B_n$ par récurrence à partir de la connaissance de $\mu$ sur les ouverts de la forme $c+p^{n-1}\Z^d$ pour $c\in B_{n-1}$. On distingue deux cas.
\begin{itemize}
\item[$\bullet$] Si $\mu\big(\mathbf{1}_{c+p^{n-1}\Z^d}(z) \prod_{i=1}^d z_i^{j_i}\big) \neq 0$ on pose pour $b\in B_n$, $b \in c+p^{n-1}\Z^d$:
\[
\mu\Big(\mathbf{1}_{b+p^n\Z^d}(z) \prod_{i=1}^d z_i^{j_i}\Big) = \left\{ \begin{array}{lll}
p^{-nr+n\sum j_i} \quad &\mbox{si} \ b=c;   \\
\mu\big(\mathbf{1}_{c+p^{n-1}\Z^d}(z) \prod_{i=1}^d z_i^{j_i}\big) - p^{-nr+n\sum j_i} \quad & \mbox{si} \ b= c+\alpha p^{n-1};\\
0 \quad &\mbox{sinon}.
\end{array} \right.
\]
\item[$\bullet$] Si $\mu\big(\mathbf{1}_{c+p^{n-1}\Z^d}(z) \prod_{i=1}^d z_i^{j_i}\big) = 0$ on pose pour $b\in B_n$, $b \in c+p^{n-1}\Z^d$:
\[
\mu\Big(\mathbf{1}_{b+p^n\Z^d}(z) \prod_{i=1}^d z_i^{j_i}\Big) = 0.
\]
\end{itemize}
Par construction la distribution $\mu$ vérifie les conditions de compatibilité \eqref{reladicomp} et la condition \eqref{primafac}. De plus on a:
\[
\big|\mu\big(\mathbf{1}_{X_n}(z) \big)\big| = \big|\mu\big(\mathbf{1}_{p^n\Z^d}(z) \big)\big| = q^{enr} = q^{enr_k} q^{en\sum_{i\neq k}r_i}
\] 
ce qui implique que la condition \eqref{secondafac} n'est pas satisfaite car $q^{enr_k}$ tend vers $+\infty$ quand $h$ tend vers $+\infty$, d'où le résultat.      
\end{proof}
\end{prop}

\begin{cor}
Les espaces de Banach $C^r(\OF,E)$ et $C^{\vec{r}}(\OF,E)$ ne sont pas isomorphes.
\begin{proof}
C'est une consequence immédiate de la Proposition \ref{noniso}.
\end{proof}
\end{cor}

\appendix\section{}

Conservons les notations du §\ref{classe} et désignons par $e_{\sigma}$ le vecteur de $\mathbb{Z}_{\geq 0}^{|S|}$ ayant toutes ses composantes nulles sauf celle d'indice $\sigma$ qui est égal à $1$. Pour $N\in \mathbb{Z}_{\geq 0}$ notons $\mathcal{O}^N(\OF,S)$ la $E$-algèbre des fonctions $P\colon \OF \to E$ définies par:
\[
P(z) = \sum_{\underline{m}\in I_{\leq N}} a_{\underline{m}} z^{\underline{m}}.
\]
Nous allons donner une estimation sur les $|a_{\underline{m}}|$ pour $\underline{m} \in I_{= N}$. 
Définissons pour $\tau \in S$ et $h \in \mathbb{Z}_{\geq 0}$ l'opérateur $\Delta_{\tau,h}\colon \mathcal{O}^N(\OF,S) \to \mathcal{O}^N(\OF,S)$ par:
\[
\Delta_{\tau,h}P (z) = \sum_{\underline{m} \in I_{\leq N}} a_{\underline{m}} \prod_{
\sigma \in S - \{\tau\}
} \sigma(z)^{m_{\sigma}} \, \tau(z+\varpi_F^h)^{m_{\tau}} - P(z).
\]
Les propriétés suivantes sont immédiates.
\begin{itemize} 
\item[$\bullet$] $\Delta_{\sigma,h} (\Delta_{\tau,h}P) = \Delta_{\tau,h}(\Delta_{\sigma,h}P)$ pour tout $\sigma,\tau \in S$ tout $h \in \mathbb{Z}_{\geq 0}$ et tout $P \in \mathcal{O}^N(\OF,S)$.
\item[$\bullet$] $\Delta_{\sigma,h} (\lambda P + \mu Q) = \lambda \Delta_{\sigma,h}P + \mu \Delta_{\sigma,h}Q$ pour tout $\sigma,\tau \in S$, tout $P,Q \in \mathcal{O}^N(\OF,S)$ et tout $\lambda, \mu \in E$.
\end{itemize}
Choisissons une numérotation $\sigma_1,\ldots, \sigma_{|S|}$ des plongements de $F$ dans $E$ et définissons pour $\underline{i} \in I_{\leq N}$ et $h \in \mathbb{Z}_{\geq 0}$ l'opérateur  $\Delta_{\underline{i},h} \colon \mathcal{O}^N(\OF,S) \to \mathcal{O}^N(\OF,S)$ par:
\[
\Delta_{\underline{i},h}P =  \Delta_{\sigma_1,h}^{i_{\sigma_1}} \circ \ldots \circ \Delta_{\sigma_{|S|},h}^{i_{\sigma_{|S|}}}P.
\]
Par la première propriété l'opérateur $\Delta_{\underline{i},h}$ ne dépend pas de la numérotation choisie.

\begin{prop}\label{ultratec}
Soit $N \in \mathbb{Z}_{\geq 0}$ et $P$ un élément de $\mathcal{O}^N(\OF,S)$ qui s'écrit sous la forme
\[
P(z) = \sum_{\underline{i}\in I_{\leq N}} a_{\underline{i}} \frac{z^{\underline{i}}}{\underline{i}!}.
\]
Alors pour tout $\underline{m} \in I_{=N}$, tout $h \in \mathbb{Z}_{\geq 0}$ et tout $z \in \OF$ on a:
\[
a_{\underline{m}} = \varpi_F^{-\underline{m}h} \Delta_{\underline{m},h}P(z).
\]
\begin{proof}
Il suffit d'étudier l'action de l'opérateur $\Delta_{\underline{m},h}$ sur les fonctions de la forme
\[
*^{\underline{i}} \colon z \mapsto z^{\underline{i}},
\]
pour ${\underline{i}} \in I_{\leq N}$. Soit $\tau \in S$. On a pour tout $z \in \OF$: 
\begin{align*}
\Delta_{e_{\tau},h}(*^{\underline{i}})(z) &= \prod_{
\sigma \in S -\{\tau\}
} \sigma(z)^{i_{\sigma}} \, \tau(z+\varpi_F^h)^{i_{\tau}} - \prod_{\sigma \in S}\sigma(z)^{i_{\sigma}} \\
&= \sum_{n=1}^{i_{\tau}} \tbinom{i_{\tau}}{n} \tau(\varpi_F^h)^n \prod_{
\sigma \in S - \{\tau\}
} \sigma(z)^{i_{\sigma}}\, \tau(z)^{i_{\tau}-n}. 
\end{align*}
En itérant la formule ci-dessus et en utilisant (ii) on déduit
\[
\forall z \in \OF, \quad \Delta_{m_{\tau}e_{\tau},h}(*^{\underline{i}})(z) =
\left\{ \begin{array}{ll}
i_{\tau}! \tau(\varpi_F^h)^{i_{\tau}} \prod_{\substack{
\sigma \in S \\
\sigma \neq \tau
}} \sigma(z)^{i_{\sigma}}   & \mbox{si} \ m_{\tau} = i_{\tau} \\
0 & \mbox{si} \ m_{\tau} > i_{\tau}.
\end{array} \right.
\]
Cela implique
\[
\forall z \in \OF, \quad \Delta_{\underline{m},h}(*^{\underline{i}})(z) =
\left\{ \begin{array}{ll}
\underline{i}! (\varpi_F^h)^{\underline{i}}   & \mbox{si} \ \underline{m} = \underline{i} \\
0 & \mbox{si} \ \underline{m} > \underline{i}.
\end{array} \right.
\]
Par ce qui précède on a pour tout $z \in \OF$:
\begin{align*}
\Delta_{\underline{m},h}P(z) = \sum_{\underline{i} \in I_{\leq N}} \frac{a_{\underline{i}}}{\underline{i}!} \Delta_{\underline{m},h}(*^{\underline{i}})(z) &= \frac{a_{\underline{m}}}{\underline{m}!} \Delta_{\underline{m},h}(*^{\underline{m}})(z) \\
&=     \frac{a_{\underline{m}}}{\underline{m}!}  \underline{m}! (\varpi_F^h)^{\underline{m}}.
\end{align*}
Cela permet de conclure.
\end{proof}
\end{prop}

\begin{lemma}\label{lemmacorollariotec}
Soit $N_1,N_2 \in \mathbb{Z}_{\geq 0}$ avec $N_2 \leq N_1$ et $P \in \mathcal{O}^{N_1}(\OF,S)$ défini par:
\[
 P(z) = \sum_{\substack{
 \underline{i}  \in I_{\leq N_1} \\
\underline{i} \in I_{\geq N_2}  
}} a_{\underline{i}} \prod_{\sigma \in S} \sigma(z)^{i_{\sigma}}.
\]
Alors il existe $C \in \mathbb{R}_{> 0}$ tel que pour tout $h \in \mathbb{Z}_{\geq 0}$ on a:
\[
\sup_{z \in \varpi_F^h \OF} |P(z)| \geq C q^{-h N_2} \sup_{z \in \OF} |P(z)|.
\]
\begin{proof}
Notons $P_k$ pour $N_2 \leq k \leq N_1$ le terme de $P$ défini par:
\[
P_k(z) = \sum_{ \underline{i} \in I_{= k}  
} a_{\underline{i}} \prod_{\sigma \in S} \sigma(z)^{i_{\sigma}},
\]
de sorte que
\[
P(z) = \sum_{k= N_2}^{N_1} P_k(z).
\]
Le cas $N_2 = N_1$ est facile et est laissé au lecteur. Supposons donc $N_1 > N_2$ et $P_{N_2} \neq 0$. Il est clair qu'il suffit de vérifier le lemme à partir d'un entier positif $h_0$ suffisamment grand. Rappelons que $e$ désigne l'indice de ramification de $F$ et qu'il existe $u \in \OF^{\times}$ tel que $u  \varpi_F^e = p$. Pour tout $h \in \mathbb{Z}_{\geq 0}$ il existe un  $\alpha_h \in \mathbb{Z}_{\geq 0}$ et $0\leq \beta_h \leq e-1$ tels que $h = \alpha_h e + \beta_h$. Posons:
\begin{align*}
M_1 &= \inf_{1 \leq l \leq e-1} \sup_{z \in \varpi_F^l \OF}|P_{N_2}(z)|, \\
M_2 & = \sup_{1 \leq l \leq e-1} \sup_{N_2+1 \leq k \leq N_1} \sup_{z \in \varpi_F^l \OF}|P_{k}(z)|.
\end{align*}
D'après les égalités:
\[
\varpi_F^h \OF = \varpi_F^{e\alpha_h} \varpi_F^{\beta_h} \OF = u^{\alpha_h}\varpi_F^{e\alpha_h} \varpi_F^{\beta_h} \OF
\]
et comme pour tout $\sigma \in S$ on a $\sigma(u \varpi_F^e) = u \varpi_F^e$, on déduit la minoration suivante:
\begin{align*}
\sup_{z \in \varpi_F^h \OF} \Big| \sum_{k=N_2 + 1}^{N_1} P_k(z)\Big| &\leq \sup_{N_2 + 1 \leq k \leq N_1} \sup_{z \in \varpi_F^h \OF} |P_k(z)| \\
 &\leq  \sup_{N_2 + 1 \leq k \leq N_1} |(u\varpi_F^{e})^{\alpha_h k}| \sup_{z \in \varpi_F^{\beta_h}\OF} |P_k(z)| \\
 &\leq  \sup_{N_2 + 1 \leq k \leq N_1} q^{-e\alpha_h k} M_2 \\
 &\leq M_2 q^{-e \alpha_h (N_2+1)}.
\end{align*}
Un calcul analogue montre que:
\[
\sup_{z \in \varpi_F^h \OF} |P_{N_2}(z)| \geq M_1 q^{-e \alpha_h N_2}. 
\]
Par ce qui précède on déduit qu'il existe un $h_0 \in \mathbb{Z}_{\geq 0}$ tel que pour tout $h \geq h_0$ on a:
\begin{align}\label{ink1}
\sup_{z \in \varpi_F^{h} \OF}|P_{N_2}(z)| > \sup_{z \in \varpi_F^{h} \OF} \Big| \sum_{k=N_2 +1}^{N_1}P_{N_k}(z)\Big|.
\end{align}
Pour tout $h \in \mathbb{Z}_{\geq 0}$ il existe $0\leq \gamma_h \leq e-1$ tel que $h+\gamma_h$ est multiple de $e$. Alors en utilisant le fait que $\sigma\Big(u^{\frac{h+\gamma_h}{e}} {\varpi_F^{h+\gamma_h}}\Big) = {u^{\frac{h+\gamma_h}{e}}\varpi_F^{h+\gamma_h}}$ pour tout $\sigma \in S$ et l'inégalité \eqref{ink1}, on a pour tout $h \geq h_0$:
\begin{align*}
\sup_{z \in \varpi_F^h \OF}|P(z)| \geq \sup_{z \in \varpi_F^{h+\gamma_h} \OF}|P(z)| &= \sup_{z \in \varpi_F^{h+\gamma_h} \OF}|P_{N_2}(z)| \\
& = q^{-(h+\gamma_h)N_2} \sup_{z \in  \OF}|P_{N_2}(z)| \\
& \geq C q^{-h N_2} \sup_{z \in  \OF}|P(z)|.
\end{align*}
où l'on a posé: 
\[
C = q^{-(e-1)N_2} \frac{\sup_{z \in \OF}|P_{N_2}(z)|}{\sup_{z \in \OF}|P(z)|}.
\]
Cela permet de conclure.  


\end{proof}
\end{lemma}

\begin{prop} \label{ultratec1}
Conservons les notations de la Proposition \ref{ultratec}. Il existe une constante $C \in \mathbb{R}_{\geq 1}$ et un $n_0 \in \mathbb{Z}_{\geq 0}$ tels que pour tout $h\geq n_0$ et tout $\underline{m} \in I_{=N}$ on a:
\[
|a_{\underline{m}}|  q^{-hN} \leq C \sup_{z \in \varpi_F^h\OF}|P(z)|.
\]

\begin{proof}
Notons $N_2$ l'entier positif tel que $P$ s'écrit sous la forme:
\[
P(z) = \sum_{\substack{
 \underline{i}  \in I_{\leq N} \\
\underline{i} \in I_{\geq N_2}  
}} a_{\underline{i}} \frac{z^{\underline{i}}}{\underline{i}!}.
\]
Montrons qu'il existe un $n_0 \in \mathbb{Z}_{\geq 0}$ et une constante $C \in \mathbb{R}_{>0}$ tels que pour tout $h \geq n_0$ et tout $\tau \in S$ on a:
\begin{align}\label{sfrb}
\sup_{z \in \varpi_F^{h}\OF} |\Delta_{e_{\tau},h}P(z)| \leq C \sup_{z \in \varpi_F^h\OF}|P(z)|.
\end{align}
Rappelons que $\Delta_{e_{\tau},h}P$ est définie pour tout $h \in \mathbb{Z}_{\geq 0}$ et tout $z \in \OF$ par:
\begin{align*}
\Delta_{e_{\tau},h}P(z) &= \sum_{\underline{i} \in I_{\leq N}} \frac{a_{\underline{i}}}{\underline{i}!} \prod_{
\sigma \in S-\{\tau\}
} \sigma(z)^{i_{\sigma}} \, \tau(z+\varpi_F^h)^{i_{\tau}} - P(z) \\
&= \sum_{\underline{i} \in I_{\leq N}} \sum_{k=1}^{i_{\tau}} \frac{a_{\underline{i}}}{\underline{i}!} \tbinom{i_{\tau}}{k} \tau(\varpi_F)^{hk} \tau(z)^{i_{\tau}-k} \prod_{\sigma \in S-\{\tau \}} \sigma(z)^{i_{\sigma}}.
\end{align*}  
On voit facilement qu'il existe $n' \in \mathbb{Z}_{\geq 0}$ tel que pour tout $\underline{i} \in I_{\leq N}$ et tout $1\leq k \leq i_{\tau}$ on a:
\begin{align}\label{disetriv}
\sup_{z \in \OF}\Big|\frac{a_{\underline{i}}}{\underline{i}!} \tau(\varpi_F)^{n'k} \tau(z)^{i_{\tau}-k} \prod_{\sigma \in S-\{\tau \}} \sigma(z)^{i_{\sigma}} \Big| \leq \sup_{z \in \OF} |P(z)|.
\end{align}
De plus, d'après le Lemme \ref{lemmacorollariotec}  on sait  qu'il existe une constante $C_1 \in \mathbb{R}_{>0}$ telle que pour tout $h \in \mathbb{Z}_{\geq 0}$ on a:
\begin{align}\label{flineg}
\sup_{z \in \varpi_F^h \OF} |P(z)| \geq C_1 q^{-h N_2} \sup_{z \in \OF} |P(z)|.
\end{align}
Notons $n''$ le plus petit entier positif tel que $|\varpi_F^{n''}| \leq C_1$. On a pour tout $h \in \mathbb{Z}_{\geq 0}$ et tout $1\leq k \leq i_{\tau}$:
\begin{align*}
&\sup_{z \in \varpi_F^{h+n'+n''}\OF}\Big|\frac{a_{\underline{i}}}{\underline{i}!} \tau(\varpi_F)^{(h+n'+n'')k} \tau(z)^{i_{\tau}-k} \prod_{\sigma \in S-\{\tau \}} \sigma(z)^{i_{\sigma}} \Big|  \\
\stackrel{\phantom{\eqref{flineg}}}=& \quad \Big|\tau(\varpi_F)^{(h+n'+n'')k} \tau \Big(\varpi_F^{(h+n'+n'')}\Big)^{i_{\tau}-k} \prod_{\sigma \in S-\{\tau\}}\sigma \Big(\varpi_F^{(h+n'+n'')}\Big)^{i_{\sigma}}\Big| \big| \tau(\varpi_F)^{-n'k} \big| \\
\cdot &   \sup_{z \in \varpi_F^{h+n'+n''}\OF}   \Big|\frac{a_{\underline{i}}}{\underline{i}!} \tau(\varpi_F)^{n'k} \tau\Big(\frac{z}{\varpi_F^{h+n'+n''}}\Big)^{i_{\tau}-k} \prod_{\sigma \in S-\{\tau\}}\sigma\Big(\frac{z}{\varpi_F^{h+n'+n''}} \Big)^{i_{\sigma}}     \Big| \\
\stackrel{\eqref{disetriv}} \leq & \quad  C C_1 q^{-N_2(h+n'+n'')} \sup_{z \in \OF} |P(z)| \\
\stackrel{\eqref{flineg}}\leq & \quad  C \sup_{z \in \varpi_F^{h+n'+n''}\OF} |P(z)|
\end{align*}
où l'on a posé $C = \sup_{1\leq k \leq m_{\tau}} \frac{1}{|\varpi_F^{(n'+n'')k}|}$. On en déduit \eqref{sfrb} pour tout $h \geq n_0$, où $n_0 \ugu n'+n''$. Donc, quitte à changer $n_0$ et $C$, on a pour tout $h\geq n_0$:
\begin{align}\label{disdif}
\sup_{z \in \varpi_F^{h}\OF} |\Delta_{\underline{m},h}P(z)| \leq C \sup_{z \in \varpi_F^h\OF}|P(z)|.
\end{align}
Or, d'après la Proposition \ref{ultratec}, on a pour tout $h \geq n_0$: 
\[
\Delta_{\underline{m},h}P(z) = a_{\underline{m}}(\varpi_F^{h})^{\underline{m}}.
\] 
et, par l'inégalité \eqref{disdif} on obtient:
\[
|a_{\underline{m}}| \, q^{-h|\underline{m}|} \leq C \sup_{z \in \varpi_F^h\OF}|P(z)|, 
\]
d'où le résultat.
\end{proof}
\end{prop}



\end{document}